\documentclass[11pt]{article}
\usepackage{amssymb}
\usepackage{epsfig}
\newtheorem{Th}{Theorem}[section]

\def\be{\begin{equation}}
\def\ee{\end{equation}}
\def\bea{\begin{eqnarray}}
\def\eea{\end{eqnarray}}
\def\bes{\begin{eqnarray*}}
\def\ees{\end{eqnarray*}}

\def\nn{\nonumber}
\def\<{\langle}
\def\>{\rangle}
\def\lb{\label}
\def\bs{\setminus}

\def\R{{\bf R}}

\def\Z{{\bf Z}}
\def\N{{\bf N}}
\def\U{{\bf U}}
\def\Q{{\bf Q}}
\def\T{{\bf T}}

\def\aa{{\alpha}}
\def\bb{{\beta}}
\def\ga{{\gamma}}

\def\th{{\theta}}
\def\Th{{\Theta}}

\def\Om{{\Omega}}
\def\ep{{\epsilon}}

\def\Lm{{\Lambda}}

\def\sg{{\sigma}}

\def\Sg{{\Sigma}}
\def\vf{{\varphi}}

\def\vth{{\vartheta}}
\def\dm{{\diamond}}

\def\H{{\cal H}}
\def\T{{\cal T}}

\def\Nn{{\cal N}}

\def\mul{{\rm mul}}

\def\crit{{\rm crit}}

\def\Sp{{\rm Sp}}

\def\per{{\rm per}}

\def\hb{\vrule height0.18cm width0.14cm $\,$}

\def\wtd#1{\widetilde{#1}}

\title{ Resonance identities and stability of symmetric closed characteristics on symmetric compact
star-shaped hypersurfaces}
\author{ Hui Liu$^{1}$,\thanks{Partially supported by NSFC (No. 11401555), China Postdoctoral Science Foundation No. 2014T70589, CUSF(No. WK3470000001).
E-mail:huiliu@ustc.edu.cn } \qquad  Yiming
Long$^{2}$,\thanks{Partially supported by NSFC (No. 11131004), MCME and LPMC of MOE of China, Nankai
University and BCMIIS of Capital Normal University. E-mail: longym@nankai.edu.cn. }  \\ \\
$^{1}$ Key Laboratory of Wu Wen-Tsun Mathematics, Chinese Academy of Sciences,\\
School of Mathematical Sciences, University of Science and Technology of China,
\\Hefei, Anhui 230026, People's Republic of China\\
$^{2}$ Chern Institute of Mathematics and LPMC, Nankai University, \\Tianjin 300071, People's Republic of China\\}
\date{}

\begin{document}

\maketitle

\begin{abstract}
{\it So far, it is still unknown whether all the closed characteristics on a symmetric compact star-shaped
hypersurface $\Sigma$ in ${\bf R}^{2n}$ are symmetric. In order to understand behaviors of such orbits,
in this paper we establish first two new resonance identities for symmetric closed characteristics on
symmetric compact star-shaped hypersurface $\Sigma$ in ${\bf R}^{2n}$ when there exist only finitely many
geometrically distinct symmetric closed characteristics on $\Sigma$, which extend the identity established
by Liu and Long in \cite{LLo1} of 2013 for symmetric strictly convex hypersurfaces. Then as an application
of these identities and the identities established by Liu, Long and Wang recently in \cite{LLW1} for all
closed characteristics on the same hypersurface, we prove that if there exist exactly two geometrically
distinct closed characteristics on a symmetric compact star-shaped hypersuface in ${\bf R}^4$, then both
of them must be elliptic.}
\end{abstract}



\renewcommand{\theequation}{\thesection.\arabic{equation}}
\renewcommand{\thefigure}{\thesection.\arabic{figure}}

\setcounter{equation}{0}%\setcounter{figure}{0}
\section{Introduction and main results}%{Section 1}

Let $\Sigma$ be a $C^3$ compact hypersurface in $\R^{2n}$ strictly star-shaped with respect to the origin,
i.e., the tangent hyperplane at any $x\in\Sigma$ does not intersect the origin. We denote the set of all
such hypersurfaces by $\H_{st}(2n)$, and denote by $\H_{con}(2n)$ the subset of $\H_{st}(2n)$ which consists
of all strictly convex hypersurfaces. We also denote the set of all hypersurfaces $\Sg\in\H_{st}(2n)$ (or
$\H_{con}(2n)$), which are symmetric with respect to the origin, i.e., $\Sigma=-\Sigma$, by
$\mathcal{S}\mathcal{H}_{st}(2n)$ (or $\mathcal{S}\mathcal{H}_{con}(2n)$). We consider closed characteristics
$(\tau, y)$ on $\Sigma$, which are solutions of the following problem
\be
\left\{\matrix{\dot{y}=JN_\Sigma(y), \cr
               y(\tau)=y(0), \cr }\right. \lb{1.1}\ee
where $J=\left(\matrix{0 &-I_n\cr
        I_n  & 0\cr}\right)$, $I_n$ is the identity matrix in $\R^n$, $\tau>0$, $N_\Sigma(y)$ is the outward
normal vector of $\Sigma$ at $y$ normalized by the condition $N_\Sigma(y)\cdot y=1$. Here $a\cdot b$ denotes
the standard inner product of $a, b\in\R^{2n}$. A closed characteristic $(\tau, y)$ is {\it prime}, if $\tau$
is the minimal period of $y$. Two closed characteristics $(\tau, y)$ and $(\sigma, z)$ are {\it geometrically
distinct}, if $y(\R)\not= z(\R)$. We denote by $\T(\Sg)$ the set of all geometrically distinct closed
characteristics on $\Sg$. A closed characteristic $(\tau, y)$ on $\Sigma \in \mathcal{S}\mathcal{H}_{st}(2n)$
is {\it symmetric} if $y({\bf R})=-y({\bf R})$, {\it non-symmetric} if $y({\bf R})\cap (-y({\bf R}))=\emptyset$.
We denote by $\mathcal{T}(\Sigma)$ ($\mathcal{T}_s(\Sigma)$) the set of geometrically distinct (symmetric)
closed characteristics $(\tau, y)$ on $\Sigma\in\mathcal {S}\mathcal{H}_{st}(2n)$. A closed characteristic
$(\tau,y)$ is {\it non-degenerate} if $1$ is a Floquet multiplier of $y$ of precisely algebraic multiplicity
$2$, and is {\it hyperbolic} if $1$ is a double Floquet multiplier of it and all the other Floquet multipliers
are not on ${\bf U}=\{z\in {\bf C}\mid |z|=1\}$, i.e., the unit circle in the complex plane, and is {\it elliptic}
if all the Floquet multipliers of $y$ are on ${\bf U}$.

The pioneer global result $\;^{\#}\mathcal{T}(\Sigma)\ge 1$ was proved by Rabinowitz in \cite{Rab1} for
$\Sg\in\mathcal{H}_{st}(2n)$ and by A. Weinstein in \cite{Wei1} for $\Sg\in\mathcal{H}_{con}(2n)$ in 1978.
For further results on the multiplicity of closed characteristics on $\Sg\in\H_{con}(2n)$ or $\H_{st}(2n)$,
we refer to \cite{EkL1}, \cite{EkH1}, \cite{Szu1}, \cite{HWZ1}, \cite{LoZ1}, \cite{WHL1}, \cite{Wan2} and
\cite{HuL1} as well as \cite{Lon3}. Recently $\;^{\#}\T(\Sg)\ge 2$ was first proved for every
$\Sg\in \H_{st}(4)$ by Cristofaro-Gardiner and Hutchings in \cite{CGH1} without any pinching or
non-degeneracy conditions. Different proofs of this result can be found in \cite{GHHM}, \cite{LLo2} and
\cite{GiG1}.

Note that for any $\Sg\in \mathcal{S}\mathcal{H}_{con}(2n)$, it was proved by Liu, Long and Zhu in
\cite{LLZ1} of 2002 that for any $(\tau,y)\in \T(\Sg)$, either $(\tau,y)$ is symmetric and then it satisfies
$y(t+\tau/2)=-y(t)$ for all $t \in \R$, or $(\tau,y)$ is non-symmetric and then $(\tau,y)\neq(\tau,-y)\in \T(\Sg)$ holds
too. Thus closed characteristics on $\Sg$ are classified into two classes, symmetric or
non-symmetric as defined above. It is natural to conjecture that $\mathcal{T}(\Sigma) = \mathcal{T}_s(\Sigma)$
for every $\Sg\in \mathcal{S}\mathcal{H}_{st}(2n)$. We are aware of only two results in this spirit. The first
one is $\;^{\#}\mathcal{T}_s(\Sigma)\ge 2$ for every $\Sg\in \mathcal{S}\mathcal{H}_{con}(2n)$ proved by Wang
(Theorem 1.1 of \cite{Wan3}, 2012). Together with \cite{LLZ1}, it implies
\be  \mathcal{T}(\Sigma) = \mathcal{T}_s(\Sigma)  \lb{1.2}\ee
for all $\Sg\in \mathcal{S}\mathcal{H}_{con}(2n)$ with $n=2$ or $3$ provided $\;^{\#}\mathcal{T}(\Sigma) = n$.
The second result is that (\ref{1.2}) holds for all $\Sg\in \mathcal{S}\mathcal{H}_{con}(8)$ provided
$\;^{\#}\mathcal{T}(\Sigma)=4$ proved by Liu, Long, Wang and Zhang recently in \cite{LLWZ}. But in general,
whether (\ref{1.2}) holds for every $\Sg\in \mathcal{S}\mathcal{H}_{st}(2n)$ is still open.

Note that recently some resonance identities for closed characteristics on compact star-shaped hypersurface
in $\R^{2n}$ were established in \cite{LLW1}, they are useful tools for studying the multiplicity and
stability of closed characteristics. Note that also recently in \cite{LLo1}, the authors established
a new resonance identity for symmetric closed characteristics on symmetric convex Hamiltonian hypersurfaces
and obtained some new results about the multiplicity and stability of symmetric orbits as applications.
Motivated by the methods of \cite{LLo1}, \cite{Vit2} and \cite{LLW1}, the first goal of this paper is to
establish two new resonance identities for symmetric closed characteristics on symmetric compact star-shaped
hypersurfaces in Theorem 1.1 below. We believe that such identities will play important roles in understanding
behaviors of closed characteristics on such symmetric hypersurfaces, including the study on (\ref{1.2}) in
general.

{\bf Theorem 1.1.} {\it Suppose $\Sigma \in \mathcal {S}\mathcal{H}_{st}(2n)$ satisfies
$\,^{\#}\mathcal{T}_s (\Sigma) < +\infty$. Denote all the geometrically distinct prime symmetric closed
characteristics on $\Sigma$ by $\{(\tau_j , y_j)\}_{1\leq j\leq k}$. Then the following identities
hold
\bea \sum_{1\le j\le k\atop \hat{\bar{i}}(y_j)>0}\frac{\hat{\bar{\chi}}(y_j)}{\hat{\bar{i}}(y_j)} &=& 1,  \lb{1.3}\\
\sum_{1\le j\le k\atop \hat{\bar{i}}(y_j)<0}\frac{\hat{\bar{\chi}}(y_j)}{\hat{\bar{i}}(y_j)} &=& 0,  \lb{1.4}\eea
where $\hat{\bar{i}}(y_j)\in\R$ is the mean index of $y_j$ given by Definition 3.6 below,
$\hat{\bar{\chi}}(y_j)\in\Q$ is the average Euler characteristic given by Definition 3.7 and
Remark 3.8 below. Specially by (3.19) below we have
\be \hat{\bar{\chi}}(y) = \frac{2}{\bar{K}(y)}\sum_{1\le k\le \bar{K}(y)/2\atop 0\le l\le 2n-2}
  (-1)^{\bar{i}(y^{2k-1})+l}\bar{k}_l(y^{2k-1}),  \lb{1.5}\ee
$\bar{K}(y)\in 2\N$ is the minimal period of critical modules of iterations of $y$ defined in Proposition
3.5, $\bar{i}(y^{m})$ is the index defined in Definition 3.6 (cf. Definition 2.7 and Theorem 2.8 below),
$\bar{k}_l(y^{m})$ is the critical
type numbers of $y^m$ given by Definition 3.2 and Remark 3.3 below. }

Note that the idea of the proof of Theorem 1.1 is a natural extension of that of Theorem 1.1 in \cite{LLW1}.
Thus in some parts of the proof of Theorem 1.1 below, we only point out its difference from and make
necessary modifications on that in \cite{LLW1}.

{\bf Remark 1.2.} When $\Sigma\in\H_{con}(2n)$ and is symmetric, we can choose $K=0$ and $\varphi$ to
satisfy Proposition 2.2 (iv) of \cite{WHL1} and (\ref{2.19}) below.
Then $e(K)=0$ in Theorem 2.16.
Noticing that $C_{S^1, l}(\bar{F}_{K}, S^1\cdot\bar{x})$ is exactly isomorphic to
$C_{S^1, l}(\Psi_a, S^1\cdot \dot{\bar{x}})$ which is defined in Definition 2.5 of \cite{LLo1}, then our identity
(\ref{1.3}) coincides with the identity (1.2) of Theorem 1.1 of \cite{LLo1}. Thus our Theorem 1.1 generalizes
the resonance identity in \cite{LLo1} for symmetric convex hypersurfaces to symmetric star-shaped hypersurfaces.

As applications of Theorem 1.1 to closed characteristics on $\Sigma \in \mathcal{S}\mathcal{H}_{st}(2n)$, we
study the stability of such orbits. For this stability problem, we refer the readers to \cite{Eke1}, \cite{DDE1},
\cite{Lon1}, \cite{Lon2}, \cite{LoZ1}, \cite{WHL1}, \cite{Wan1} and the references therein. In particular, in
\cite{Lon2} of 2000, Long proved that $\Sg\in\H_{con}(4)$ and $\,^{\#}\T(\Sg)=2$ imply that both of the closed
characteristics must be elliptic, in \cite{WHL1} of 2007, W. Wang, X. Hu and Y. Long proved further that
$\Sg\in\H_{con}(4)$ and $\,^{\#}\T(\Sg)=2$ imply that both of the closed characteristics must be irrationally
elliptic, i.e., each of them possesses four Floquet multipliers with two $1$'s and the other two locate on the
unit circle with rotation angles being irrational multiples of $\pi$.

{\bf Remark 1.3.} Using Theorem 1.1, and the proof of Theorem 1.4 of \cite{LLo1} or the proof of Theorem 1.1 of
\cite{Vit2}, we obtain the following immediately:

{\it In the $C^\infty$ topology, the following holds for a generic
$\Sigma\in\mathcal {S}\mathcal {H}_{st}(2n)$: the problem (1.1) has infinitely many prime symmetric
closed characteristics, or all the symmetric closed characteristics on $\Sigma$ are hyperbolic.}

Note that for a $\Sigma \in \mathcal{S}\mathcal{H}_{st}(2n)\bs \mathcal{S}\mathcal{H}_{con}(2n)$, we are not
aware of any results concerning the existence of at least two elliptic closed characteristics without pinching
or non-degeneracy condition. Thus motivated by result of \cite{Lon2}, we prove the following result as an
application of Theorem 1.1.

{\bf Theorem 1.4.} {\it Let $\Sg\in \mathcal {S}\mathcal{H}_{st}(4)$ satisfy $\,^{\#}\T(\Sg)=2$. Then both of
the closed characteristics must be elliptic.}

{\bf Remark 1.5.} Note that the symmetric condition on $\Sg$ in Theorem 1.4 is specially used to get a
contradiction in the study of the Subcase 1.2 of Case 1 in the proof of Theorem 1.4, where we used Theorem 1.1
to get the identity (\ref{6.36}) and then to get the estimate (\ref{6.38}) below.

Besides the resonance identities established in \cite{LLW1} and our Theorem 1.1, the other main ingredients in
the proof of Theorem 1.4 are: Morse inequality, and the index iteration theory developed by Long and his
coworkers, specially the precise iteration formulae of the Maslov-type index theory for any symplectic path
which is established by Long in \cite{Lon2} and the common index jump theorem of Long and Zhu (Theorem 4.3 of
\cite{LoZ1}).

This paper is arranged as follows. In Sections 2-4, we give first a proof for Theorem 1.1, and then in Section 5
we briefly review the equivariant Morse theory and the resonance identities for closed characteristics on
compact star-shaped hypersurfaces in $\R^{2n}$ developed in \cite{LLW1}. The proof of Theorem 1.4 will be given
in Section 6.

In this paper, let $\N$, $\N_0$, $\Z$, $\Q$, $\R$, and $\R^+$ denote the sets of natural integers,
non-negative integers, integers, rational numbers, real numbers, and positive real numbers respectively.
We define the function $[a]=\max{\{k\in {\bf Z}\mid k\leq a\}}$, $\{a\}=a-[a]$ , and
$E(a)=\min{\{k\in{\bf Z}\mid k\geq a\}}$.
Denote by $a\cdot b$ and $|a|$ the standard
inner product and norm in $\R^{2n}$. Denote by $\langle\cdot,\cdot\rangle$ and $\|\cdot\|$
the standard $L^2$ inner product and $L^2$ norm. For an $S^1$-space $X$, we denote by
$X_{S^1}$ the homotopy quotient of $X$ by $S^1$, i.e., $X_{S^1}=S^\infty\times_{S^1}X$,
where $S^\infty$ is the unit sphere in an infinite dimensional {\it complex} Hilbert space.
In this paper we use $\Q$ coefficients for all homological and cohomological modules. By $t\to a^+$, we
mean $t>a$ and $t\to a$.

\setcounter{equation}{0}
\section{Critical point theory for symmetric closed characteristics}%{Section 2}
Firstly, we briefly review the variational construction of closed characteristics on starshaped hypersurfaces established in \cite{LLW1}.
Now we fix a $\Sg\in\H_{st}(2n)$ and assume the following condition on $\T(\Sg)$:

\noindent (F) {\bf There exist only finitely many geometrically distinct prime closed characteristics
\\$\quad \{(\sigma_j, z_j)\}_{1\le j\le k^\prime}$ on $\Sigma$. }

Let $\hat{\sigma}=\inf_{1\leq j\leq k^\prime}{\sigma_j}$ and $T$ be a fixed positive constant.
Then by Section 2 of \cite{LLW1}, for any
$a>\frac{\hat{\sigma}}{T}$, we can construct a  function $\varphi_a\in C^{\infty}({\bf R}, {\bf R}^+)$ which has 0
as its unique critical point in $[0, +\infty)$. Moreover,
$\frac{\varphi^{\prime}(t)}{t}$ is strictly decreasing for $t>0$ together with $\varphi(0)=0=\varphi^{\prime}(0)$ and
$\varphi^{\prime\prime}(0)=1=\lim_{t\rightarrow
0^+}\frac{\varphi^{\prime}(t)}{t}$. More precisely, we
define $\varphi_a$ and the Hamiltonian function $\wtd{H}_a(x)=a\vf_a(j(x))$ via Lemma 2.2 and Lemma 2.4
in \cite{LLW1}. The precise dependence of $\varphi_a$ on $a$ is explained in Remark 2.3 of \cite{LLW1}.

For technical reasons we want to further modify the Hamiltonian, we define the new Hamiltonian
function $H_a$ via Proposition 2.5 of \cite{LLW1} and consider the fixed period problem
\be   \left\{\matrix{\dot{x}(t) &=& JH_a^\prime(x(t)), \cr
                   x(0) &=& x(T).\qquad         \cr }\right.  \lb{2.1}\ee
Then $H_a\in C^{3}({\bf R}^{2n} \setminus\{0\},{\bf R})\cap C^{1}({\bf R}^{2n},{\bf R})$.
Solutions of (\ref{2.1}) are $x\equiv 0$ and $x=\rho z(\sigma t/T)$ with
$\frac{\vf_a^\prime(\rho)}{\rho}=\frac{\sigma}{aT}$, where $(\sigma, z)$ is a solution of (\ref{1.1}). In particular,
non-zero solutions of (\ref{2.1}) are in one to one correspondence with solutions of (\ref{1.1}) with period
$\sigma<aT$.

For any $a>\frac{\hat{\sigma}}{T}$, we can choose some
large constant $K=K(a)$ such that
\be H_{a,K}(x) = H_a(x)+\frac{1}{2}K|x|^2   \lb{2.2}\ee
is a strictly convex function, that is,
\be (\nabla H_{a, K}(x)-\nabla H_{a, K}(y), x-y) \geq \frac{\ep}{2}|x-y|^2,  \lb{2.3}\ee
for all $x, y\in {\bf R}^{2n}$, and some positive $\ep$. Let $H_{a,K}^*$ be the Fenchel dual of $H_{a,K}$
defined by
\bea  H_{a,K}^\ast (y) = \sup\{x\cdot y-H_{a,K}(x)\;|\; x\in \R^{2n}\}.   \nn\eea
The dual action functional on $X=W^{1, 2}({\bf R}/{T {\bf Z}}, {\bf R}^{2n})$ is defined by
\be F_{a,K}(x) = \int_0^T{\left[\frac{1}{2}(J\dot{x}-K x,x)+H_{a,K}^*(-J\dot{x}+K x)\right]dt}.  \lb{2.4}\ee
Then $F_{a,K}\in C^{1,1}(X, \R)$ and for $KT\not\in 2\pi{\bf Z}$, $F_{a,K}$ satisfies the
Palais-Smale condition and $x$ is a
critical point of $F_{a, K}$ if and only if it is a solution of (\ref{2.1}). Moreover,
$F_{a, K}(x_a)<0$ and it is independent of $K$ for every critical point $x_a\neq 0$ of $F_{a, K}$.

When $KT\notin 2\pi{\bf Z}$, the map $x\mapsto -J\dot{x}+Kx$ is a Hilbert space isomorphism between
$X=W^{1, 2}({\bf R}/{T {\bf Z}}; {\bf R}^{2n})$ and $E=L^{2}({\bf R}/(T {\bf Z}),{\bf R}^{2n})$. We denote its inverse
by $M_K$ and the functional
\be \Psi_{a,K}(u)=\int_0^T{\left[-\frac{1}{2}(M_{K}u, u)+H_{a,K}^*(u)\right]dt}, \qquad \forall\,u\in E. \lb{2.5}\ee
Then $x\in X$ is a critical point of $F_{a,K}$ if and only if $u=-J\dot{x}+Kx$ is a critical point of $\Psi_{a, K}$.

Suppose $u$ is a nonzero critical point of $\Psi_{a, K}$.
Then the formal Hessian of $\Psi_{a, K}$ at $u$ is defined by
\be Q_{a,K}(v)=\int_0^T(-M_K v\cdot v+H_{a,K}^{*\prime\prime}(u)v\cdot v)dt,  \lb{2.6}\ee
which defines an orthogonal splitting $E=E_-\oplus E_0\oplus E_+$ of $E$ into negative, zero and positive subspaces.
The index and nullity of $u$ are defined by $i_K(u)=\dim E_-$ and $\nu_K(u)=\dim E_0$ respectively.
Similarly, we define the index and nullity of $x=M_Ku$ for $F_{a, K}$, we denote them by $i_K(x)$ and
$\nu_K(x)$. Then we have
\be  i_K(u)=i_K(x),\quad \nu_K(u)=\nu_K(x),  \lb{2.7}\ee
which follow from the definitions (\ref{2.4}) and (\ref{2.5}). The following important formula was proved in
Lemma 6.4 of \cite{Vit2}:
\be  i_K(x) = 2n([KT/{2\pi}]+1)+i^v(x) \equiv d(K)+i^v(x),   \lb{2.8}\ee
where the index $i^v(x)$ does not depend on K, but only on $H_a$. By Theorem 6.1 below,
we know $i^v(x)=i(x)-n$, where $i(x)$ is the Maslov-type index of the fundamental solution of the system (\ref{2.9})
below(cf. Section 5.4 and Chapter 8 of \cite{Lon3}).

By the proof of Proposition 2 of \cite{Vit1}, we have that $v\in E$ belongs to the null space of $Q_{a, K}$
if and only if $z=M_K v$ is a solution of the linearized system
\be  \dot{z}(t) = JH_a''(x(t))z(t).  \lb{2.9}\ee
Thus the nullity in (\ref{2.7}) is independent of $K$, which we denote by $\nu^v(x)\equiv \nu_K(u)= \nu_K(x)$.

By Proposition 2.11 of \cite{LLW1}, the index $i^v(x)$ and nullity $\nu^v(x)$ coincide with those defined for
the Hamiltonian $H(x)=j(x)^\alpha$ for all $x\in\R^{2n}$ and some $\aa\in (1,2)$. Especially
$1\le \nu^v(x_b)\le 2n-1$ always holds.

We have a natural $S^1$-action on $X$ or $E$ defined by
\be  \theta\cdot u(t)=u(\theta+t),\quad\forall\, \theta\in S^1, \, t\in\R.  \lb{2.10}\ee
Clearly both of $F_{a, K}$ and $\Psi_{a, K}$ are $S^1$-invariant. For any $\kappa\in\R$, we denote by
\bea
\Lambda_{a, K}^\kappa &=& \{u\in L^{2}({\bf R}/{T {\bf Z}}; {\bf R}^{2n})\;|\;\Psi_{a,K}(u)\le\kappa\},  \lb{2.11}\\
X_{a, K}^\kappa &=& \{x\in W^{1, 2}({\bf R}/(T {\bf Z}),{\bf R}^{2n})\;|\;F_{a, K}(x)\le\kappa\}.  \lb{2.12}\eea
For a critical point $u$ of $\Psi_{a, K}$ and the corresponding $x=M_K u$ of $F_{a, K}$, let
\bea
\Lm_{a,K}(u) &=& \Lm_{a,K}^{\Psi_{a, K}(u)}
   = \{w\in L^{2}(\R/(T\Z), \R^{2n}) \;|\; \Psi_{a, K}(w)\le\Psi_{a,K}(u)\},  \lb{2.13}\\
X_{a,K}(x) &=& X_{a,K}^{F_{a,K}(x)} = \{y\in W^{1, 2}(\R/(T\Z), \R^{2n}) \;|\; F_{a,K}(y)\le F_{a,K}(x)\}. \lb{2.14}\eea
Clearly, both sets are $S^1$-invariant. Denote by $\crit(\Psi_{a, K})$ the set of critical points of $\Psi_{a, K}$.
Because $\Psi_{a,K}$ is $S^1$-invariant, $S^1\cdot u$ becomes a critical orbit if $u\in \crit(\Psi_{a, K})$.

In the following, we construct a variational structure of symmetric closed characteristics.
In sections 2-4, we fix first a $\Sg\in\mathcal {S}\H_{st}(2n)$ and assume the following condition on $\T_s(\Sg)$:

\noindent $(F^\prime)$ {\bf There exist only finitely many geometrically distinct prime symmetric closed characteristics
$\{(\tau_j, y_j)\}_{1\le j\le k}$ on $\Sigma$. }

Note that $(\tau, y)$ is some odd iteration of a prime symmetric closed characteristic if and only if it satisfies
the problem
\be  \left\{\matrix{
      \dot{y}(t)=JN_{\Sigma}(y(t)),  y(t)\in \Sg, \forall\;t\in {\bf R}, \cr
              y(\tau/2)=-y(0).       \cr}\right.  \lb{2.15}\ee
As Definition 2.1 of \cite{LLo1} we introduce the following discrete subset of $\R^+$:

{\bf Definition 2.1.} {\it Under the assumption (F), the set of periods of symmetric closed characteristics on
$\Sigma$ is defined by}
$$  \per(\Sigma)=\{(2m-1)\tau_j \mid m\in {\bf N}, 1\leq j \leq k\}.   $$

Note that in the above definition, the period set $\per(\Sigma)$ is defined only via odd iterations of
$(\tau_j,y_j)$, because even iterate $(2m\tau_j,y_j)$ does not satisfy the equation (\ref{2.15}) and
does not yield any critical point of $\bar{F}_{a, K}$ on the space $\bar{X}$ via Proposition 2.4 and Lemma 2.5.

Now we construct a variational structure of symmetric closed characteristics as the following.

{\bf Lemma 2.2.}(cf. Lemma 2.2 of \cite{LLW1}) {\it For any sufficiently small $\vth\in (0,1)$, there exists a function
$\vf\equiv \vf_{\vth}\in C^\infty(\R, \R^+)$ depending on $\vth$ which has $0$ as its
unique critical point in $[0, +\infty)$ such that the following hold.

(i) $\vf(0)=0=\vf^\prime(0)$, and $\vf^{\prime\prime}(0)=1=\lim_{t\rightarrow 0^+}\frac{\vf^\prime(t)}{t}$;

(ii) $\frac{d}{dt}\left(\frac{\vf^\prime(t)}{t}\right)<0$ for $t>0$, and
$\lim_{t\rightarrow +\infty}\frac{\vf^\prime(t)}{t}<\vth$; that is, $\frac{\vf^\prime(t)}{t}$ is
strictly decreasing for $t> 0$;

(iii) In particular, we can choose $\aa\in (1,2)$ sufficiently close to $2$ and $c\in (0,1)$ such that
$\vf(t)=ct^{\aa}$ whenever $\frac{\vf'(t)}{t}\in [\vth, 1-\vth]$ and $t>0$.}

Let $j: \R^{2n}\rightarrow\R$ be the gauge function of $\Sigma$, i.e., $j(\lambda x)=\lambda$ for $x\in\Sigma$
and $\lambda\ge0$, then $j\in C^3(\R^{2n}\bs\{0\}, \R)\cap C^0(\R^{2n}, \R)$ and $\Sigma=j^{-1}(1)$.
Denote by $\hat{\tau}=\inf\{s\,|\, s\in \per(\Sigma)\}$.

By the same proofs of Lemma 2.4 and Proposition 2.5 of \cite{LLW1}, we have

{\bf Lemma 2.3.} {\it Let $a>\frac{\hat{\tau}}{T}$, $\vth_a\in (0, \frac{\hat{\tau}}{aT})$ and $\vf_a$ be a
$C^\infty$ function associated to $\vth_a$ satisfying (i)-(ii) of Lemma 2.2. Define the Hamiltonian function
$\wtd{H}_a(x)=a\vf_a(j(x))$ and consider the fixed period system
\be \left\{\matrix{\dot{x}(t) &=& J\wtd{H}_a^\prime(x(t)), \cr
                   x(T/2) &=& -x(0).\qquad    \cr }\right.   \lb{2.16}\ee
Then solutions of (\ref{2.16}) are $x\equiv 0$ and $x=\rho y(\tau t/T)$ with
$\frac{\vf_a^\prime(\rho)}{\rho}=\frac{\tau}{aT}$, where $(\tau, y)$ is a solution of (\ref{2.15}). In particular,
non-zero solutions of (\ref{2.16}) are in one to one correspondence with solutions of (\ref{2.15}) with period
$\tau<aT$.}

{\bf Proposition 2.4.} {\it For $a>\frac{\hat{\tau}}{T}$ and small $\ep_a$, we choose small enough $\vth_a$ such
that Lemma 2.3 holds. Then there exists a function $\bar{H}_a$ on ${\bf R}^{2n}$ such that
$\bar{H}_a(x)=\bar{H}_a(-x)$ for all $x\in\R^{2n}$, $\bar{H}_a$ is $C^1$ on ${\bf R}^{2n}$,
and $C^3$ on ${\bf R}^{2n}\bs\{0\}$, $\bar{H}_a=\wtd{H}_a$ in $U_A\equiv\{x\mid \wtd{H}_a(x)\leq A\}$
for some large $A$, and $\bar{H}_a(x)=\frac{1}{2}\ep_a|x|^2$ for $|x|$ large,
and the solutions of the fixed period system
\be   \left\{\matrix{\dot{x}(t) &=& J\bar{H}_a^\prime(x(t)), \cr
                   x(T/2) &=& -x(0),\qquad         \cr }\right.  \lb{2.17}\ee
are the same with those of (\ref{2.16}).}

As in \cite{BLMR} (cf. Section 3 of \cite{Vit2}), for any $a>\frac{\hat{\tau}}{T}$, we can choose some
large constant $K=K(a)$ such that
\be \bar{H}_{a,K}(x) = \bar{H}_a(x)+\frac{1}{2}K|x|^2   \lb{2.18}\ee
is a strictly convex function, that is,
\be (\nabla \bar{H}_{a, K}(x)-\nabla \bar{H}_{a, K}(y), x-y) \geq \frac{\ep}{2}|x-y|^2,  \lb{2.19}\ee
for all $x, y\in {\bf R}^{2n}$, and some positive $\ep$. Let $\bar{H}_{a,K}^*$ be the Fenchel dual of
$\bar{H}_{a,K}$ defined by
\be  \bar{H}_{a,K}^\ast (y) = \sup\{x\cdot y-\bar{H}_{a,K}(x)\;|\; x\in \R^{2n}\}.   \lb{2.20}\ee
Since $\bar{H}_a(x)=\bar{H}_a(-x)$ for all $x\in\R^{2n}$, then $\bar{H}_{a,K}(x)=\bar{H}_{a,K}(-x)$ and
$\bar{H}_{a,K}^\ast(x)=\bar{H}_{a,K}^\ast (-x)$ for all $x\in\R^{2n}$ by (\ref{2.18}) and (\ref{2.20}) respectively.

The dual action functional on $\bar{X}=W_s^{1, 2}({\bf R}/{T {\bf Z}}, {\bf R}^{2n})=\{x\in
W^{1, 2}({\bf R}/{T {\bf Z}}, {\bf R}^{2n})\mid x(t+T/2)=-x(t),\forall t\in\R\}$ is defined by
\be \bar{F}_{a,K}(x) = \int_0^{T/2}{\left[\frac{1}{2}(J\dot{x}-K x,x)+\bar{H}_{a,K}^*(-J\dot{x}+K x)
\right]dt}.\lb{2.21}\ee
Then $\bar{F}_{a,K}\in C^{1,1}(\bar{X},{\bf R})$ holds by the same argument in the proof of (3.16)
of \cite{Vit2}, but $\bar{F}_{a,K}$ is not $C^2$.

{\bf Lemma 2.5.} {\it Assume $\frac{KT}{2\pi}\not\in 2{\bf Z}-1$, then $x$ is a
critical point of $\bar{F}_{a, K}$ if and only if it is a solution of (\ref{2.17}).}

{\bf Proof.} Noticing that when $\frac{KT}{2\pi}\not\in 2{\bf Z}-1$, the map $x\mapsto -J\dot{x}+Kx$
is a Hilbert space isomorphism between
$\bar{X}=W_s^{1, 2}({\bf R}/{T {\bf Z}}; {\bf R}^{2n})$ and
$\bar{E}=L_s^{2}({\bf R}/(T {\bf Z}),{\bf R}^{2n})=\{u\in L^{2}({\bf R}/(T {\bf Z}),{\bf R}^{2n})\mid
u(t+T/2)=-u(t), a.e. t\in\R\}$. Then the lemma follows by direct computation. cf. Proposition 3.4 of \cite{Vit2}.
\hfill\hb

From Lemma 2.5, we know that the critical points of $\bar{F}_{a,K}$ are independent of $K$.

{\bf Proposition 2.6.} {\it For every critical point $x_a\neq 0$ of $\bar{F}_{a,K}$, the critical value
$\bar{F}_{a,K}(x_a)<0$ holds and is independent of $K$.}

{\bf Proof.} Since $\nabla \bar{H}_{a,K}(x_a)=-J\dot{x}_a+Kx_a$, then we have
$$  \bar{H}_{a,K}^*(-J\dot{x}_a+Kx_a)=(-J\dot{x}_a+Kx_a, x_a)-\bar{H}_{a,K}(x_a).  $$
Thus we obtain
\bea \bar{F}_{a,K}(x_a)
&=& \int_0^{\frac{T}{2}}{\left[\frac{1}{2}(J\dot{x}_a-Kx_a,x_a)+\bar{H}_{a,K}^*(-J\dot{x}_a+Kx_a)\right]dt}  \nn\\
&=& \int_0^{\frac{T}{2}}{\left[-\frac{1}{2}(J\dot{x}_a-Kx_a,x_a)-\bar{H}_{a,K}(x_a)\right]dt}  \nn\\
&=& \int_0^{\frac{T}{2}}{\left[-\frac{1}{2}(J\dot{x}_a,x_a)-\bar{H}(x_a)\right]dt}  \nn\\
&=& \int_0^{\frac{T}{2}}{\left[\frac{1}{2}(\bar{H}_a^\prime(x_a),x_a)-\bar{H}(x_a)\right]dt}.  \lb{2.22}\eea
By Lemma 2.3 and Proposition 2.4, we have $x_a=\rho_a y(\tau t/T)$ with
$\frac{\vf_a^\prime(\rho_a)}{\rho_a}=\frac{\tau}{aT}$. Hence, we have
\be  \bar{F}_{a,K}(x_a)=\frac{1}{4}a\vf_a^\prime(\rho_a)\rho_aT-\frac{1}{2}a\vf_a(\rho_a)T.   \lb{2.23}\ee
Here we used the facts that $j^\prime(y)=N_\Sigma(y)$ and $j^\prime(y)\cdot y=1$.

Let $f(t)=\frac{1}{4}a\vf_a^\prime(t)t-\frac{1}{2}a\vf_a(t)$ for $t\ge 0$. Then we have $f(0)=0$ and
$f'(t)=\frac{a}{4}(\vf_a^{\prime\prime}(t)t-\vf_a^\prime(t))<0$ since $\frac{d}{dt}(\frac{\vf_a^\prime(t)}{t})<0$
by (ii) of Lemma 2.2. Together with (\ref{2.23}), it yields the proposition. \hfill\hb

We know that when $\frac{KT}{2\pi}\not\in 2{\bf Z}-1$, the map $x\mapsto -J\dot{x}+Kx$ is a Hilbert space
isomorphism between $\bar{X}=W_s^{1, 2}({\bf R}/{T {\bf Z}}; {\bf R}^{2n})$ and
$\bar{E}=L_s^{2}({\bf R}/(T {\bf Z}),{\bf R}^{2n})$. We denote its inverse by $\bar{M}_K$ and the functional
\be \bar{\Psi}_{a,K}(u)=\int_0^{\frac{T}{2}}{\left[-\frac{1}{2}(\bar{M}_{K}u, u)+
\bar{H}_{a,K}^*(u)\right]dt}, \qquad \forall\,u\in \bar{E}. \lb{2.24}\ee
Then $x\in \bar{X}$ is a critical point of $\bar{F}_{a,K}$ if and only if $u=-J\dot{x}+Kx$
is a critical point of $\bar{\Psi}_{a, K}$.
We have a natural $S^1$-action on $\bar{X}$ or $\bar{E}$ defined by
\be  \theta\cdot x(t)=x(\theta+t),\quad\forall\, \theta\in S^1, \, t\in\R.  \lb{2.25}\ee
Then we have the functionals $\bar{F}_{a, K}$, $\bar{\Psi}_{a, K}$ are $S^1$-invariant and are
even by the eveness of $\bar{H}_{a,K}^\ast(x)$ on $\bar{X}$. cf. Lemma 2.3 of \cite{Wan3}.

For any $\kappa\in\R$, we denote by
\bea
\bar{\Lambda}_{a, K}^\kappa &=& \{u\in \bar{E}\;|\;\bar{\Psi}_{a,K}(u)\le\kappa\},  \lb{2.26}\\
\bar{X}_{a, K}^\kappa &=& \{x\in \bar{X}\;|\;\bar{F}_{a, K}(x)\le\kappa\},  \lb{2.27}\eea
which are $S^1$-invariant.

{\bf Definition 2.7.} {\it Suppose $u$ is a nonzero critical point of $\bar{\Psi}_{a, K}$.
Then the formal Hessian of $\bar{\Psi}_{a, K}$ at $u$ is defined by\be
\bar{Q}_{a,K}(v)=\int_0^{\frac{T}{2}}(-\bar{M}_K v\cdot v+\bar{H}_{a,K}^{*\prime\prime}(u)v\cdot v)dt,\lb{2.28}\ee
which defines an orthogonal splitting $\bar{E}=\bar{E}_-\oplus \bar{E}_0\oplus \bar{E}_+$ of
$\bar{E}$ into negative, zero and positive subspaces.
The index and nullity of $u$ are defined by $\bar{i}_K(u)=\dim \bar{E}_-$ and
$\bar{\nu}_K(u)=\dim \bar{E}_0$ respectively. }

Similarly, we define the index and nullity of $x=\bar{M}_Ku$ for $\bar{F}_{a, K}$ respectively by
\be  \bar{i}_K(x) = \bar{i}_K(u),\quad \bar{\nu}_K(x) = \bar{\nu}_K(u),  \lb{2.29}\ee
which follow from the definitions (\ref{2.21}) and (\ref{2.24}).

By the same proof of Proposition 2 of \cite{Vit1}, we have that $v\in \bar{E}$ belongs to the null space
of $\bar{Q}_{a, K}$ if and only if $z=\bar{M}_K v$ is a solution of the linearized system
\be  \left\{\matrix{\dot{z}(t)=J\bar{H}_{a}^{\prime\prime}(x(t))z(t), \cr
                    z(\frac{T}{2})=-z(0).   \cr }\lb{2.30}\right.\ee
Thus the nullity in (\ref{2.28}) is independent of $K$, which we denote by
$\bar{\nu}(x)\equiv \bar{\nu}_K(u)= \bar{\nu}_K(x)$.

Suppose $u$ is a nonzero critical point of $\bar{\Psi}_{a, K}$ such that $u$ corresponds to a
critical point $x=\bar{M}_Ku$ of $\bar{F}_{a, K}$. Then
for any $\omega=-1,1$ and $k=1,2$, we define the quadratic form
\bea Q_{\frac{kT}{2}}^\omega(v)=\int_0^{\frac{kT}{2}}(-\bar{M}_K v\cdot v+\bar{H}_{a,K}^{*\prime\prime}(u)v\cdot v)dt
\nn\eea on $E_{\frac{kT}{2}}^\omega=\{u\in L^{2}({\bf R}/(T {\bf Z}),{\bf R}^{2n})\mid
u(t+\frac{kT}{2})=\omega u(t), a.e. t\in\R\}$,
we denote its index by $i_{K, \frac{kT}{2}}^\omega(u)$. Then we have
\bea \bar{i}_K(u)=i_{K, \frac{T}{2}}^{-1}(u).\lb{2.31}\eea
By the same proof of Corollary 1.5.4 of \cite{Eke2}, we obtain\bea
i_{K, T}^1(u)=i_{K, \frac{T}{2}}^1(u)+i_{K, \frac{T}{2}}^{-1}(u).\lb{2.32}\eea
Consider the linear Hamiltonian system
\be \left\{\matrix{\dot{\xi}(t)=J A_a(t)\xi, \qquad {\rm for}\;\;0\leq t\leq T/2,\cr
                   A_a(T/2)=A_a(0),   \cr}\lb{2.33}\right.\ee
where $A_a(t)=\bar{H}_a^{\prime\prime}(x(t))$. Let $\psi_x:[0,T/2]\to\Sp(2n)$ with $\psi_x(0)=I_{2n}$ be the fundamental solution
of (\ref{2.33}). Denote by $i(A_a, k)\equiv i(\psi_x, k)$
and $\nu(A_a, k)\equiv \nu(\psi_x, k)$ the Maslov-type index and nullity of the $k$-th iteration of the system (\ref{2.33})
(cf. Section 5.4 and Chapter 8 of \cite{Lon3}).

Then from Lemma 6.4 of \cite{Vit2} and Theorem 2.1 of \cite{HuL1}, we have
\bea i_{K, \frac{kT}{2}}^1(u)=2n([K\frac{kT}{4\pi}]+1)+i(A_a, k)-n.\lb{2.34}\eea
In the following we set $e(K)=2n([K\frac{T}{2\pi}]-[K\frac{T}{4\pi}])$.

{\bf Theorem 2.8.} {\it Suppose $u$ is a nonzero critical point of $\bar{\Psi}_{a, K}$ such that
$u$ corresponds to a critical point $x=\bar{M}_Ku$ of $\bar{F}_{a, K}$. Then we have
\be \bar{i}_K(u)=\bar{i}_K(x)=e(K)+i_{-1}(A_a,1), \;\;\bar{\nu}_K(u)=\bar{\nu}_K(x)=\nu_{-1}(A_a,1), \lb{2.35}\ee
where $i_{-1}(A_a,1)\equiv i_{-1}(\psi_x)$ and $\nu_{-1}(A_a,1)\equiv \nu_{-1}(\psi_x)$ are the Maslov-type index and nullity introduced in Definition 5.4.3 of
\cite{Lon3} for $\omega=-1$. They depend only on $a$, and are independent of $K$. In the following, we also
denote them by $i_{-1}(x)$ and $\nu_{-1}(x)$ respectively.}

{\bf Proof.} From (\ref{2.29}), (\ref{2.31}), (\ref{2.32}), (\ref{2.34}), it follows that\bea
\bar{i}_K(u)=\bar{i}_K(x)=e(K)+i(A_a,2)-i(A_a,1).\lb{2.36}\eea
By the Bott-type formulae (cf. Theorem 9.2.1 of \cite{Lon3}), we have\bea
i(A_a,2)=i(A_a,1)+i_{-1}(A_a,1).\lb{2.37}\eea
Thus the former in (\ref{2.35}) holds by (\ref{2.36}) and (\ref{2.37}). The latter in (\ref{2.35})
follows from (\ref{2.30}).\hfill\hb

{\bf Lemma 2.9.} {\it There is an $S^1$-invariant subspace $V$ of
$\bar{X}=W_s^{1, 2}({\bf R}/{T {\bf Z}}, {\bf R}^{2n})$
such that $\bar{F}_{a,K}$ is bounded from below on $V$ and $\bar{F}_{a,K}(z)$ goes to minus infinity
as $\|z\|$ goes to infinity
on $V^\bot$, where $codim V=e(K)$, $e(K)=2n([K\frac{T}{2\pi}]-[K\frac{T}{4\pi}])$.}

{\bf Proof.} Define the quadratic form $A$ by
\bea A(z) = \int_0^{T/2}{\left[(J\dot{z}-K z,z)+\frac{1}{K+\epsilon}|-J\dot{z}+K z|^2\right]dt}.  \nn\eea
and take for $V^\bot$ the space generated by the eigenvectors associated with negative eigenvalues.

Let $z(t)=\sum_{-\infty}^{+\infty}{\exp{(\frac{2\pi}{T}(2k-1)Jt)}z_k}$ be the Fourier decomposition of $z$.
Then
\bea A(z)=\sum_{-\infty}^{+\infty}{\frac{((4k-2)\pi -\epsilon T)((4k-2)\pi+KT)}{2(K+\epsilon)T^2}|z_k|^2}, \lb{2.38}\eea
and the codimension of $V$ is given by $2n\,^{\#}\{k\in\Z\mid -KT<(4k-2)\pi<\epsilon T\} = e(K)$.
Hence by the same proof of Lemma 5.2 of \cite{Vit2}, our lemma follows.\hfill\hb

Because $\bar{X}$ is $S^1$-equivariantly homotopic to the single point $0$ in itself, as in the proof of the
identity (5.3) of \cite{LLo1}, by Lemma 2.9 and the method of Corollary 5.11 of \cite{Vit2} we obtain the
following result on the global equivariant homological structure of $(\bar{X}, \bar{X}^{-\infty})$, which
will be used in the proof of the identity (\ref{4.5}) below.

{\bf Corollary 2.10.} {\it $H_{S^1}^q(\bar{X}, \bar{X}^{-\infty})\simeq H^{q-e(K)}(CP^{\infty})$
(which is one dimensional for all even $q\geq e(K)$, and $\{0\}$ otherwise).}

In this paper, we say that $\bar{\Psi}_{a,K}$ with $a\in [a_1, a_2]$ form a continuous family of functionals
in the sense of Remark 2.6 of \cite{LLW1}, when $0<a_1<a_2<+\infty$.

{\bf Lemma 2.11.} {\it For any $0<a_1<a_2<+\infty$, let $K$ be fixed so that $\bar{\Psi}_{a, K}$ with $a\in [a_1, a_2]$
is a continuous family of functionals defined by (\ref{2.24}) satisfying (\ref{2.19}) with the same $\ep>0$. Then
there exist a finite dimensional $S^1$-invariant subspace $\bar{G}$ of
$\bar{E}=L_s^{2}({\bf R}/(T {\bf Z}),{\bf R}^{2n})$ and
a family of $S^1$-equivariant maps $\bar{h}_{a}: \bar{G}\rightarrow \bar{G}^\perp$ such that the following hold.

(i) For $g\in \bar{G}$, each function $h\mapsto\bar{\Psi}_{a,K}(g+h)$ has $\bar{h}_a(g)$ as
the unique minimum in $\bar{G}^\perp$.

Let $\bar{\psi}_{a,K}(g)=\bar{\Psi}_{a,K}(g+\bar{h}_a(g))$. Then we have

(ii) Each $\bar{\psi}_{a, K}$ is $C^1$ and $S^1$-invariant on $\bar{G}$. Here $g_a$ is a critical point
of $\bar{\psi}_{a, K}$ if and only if $g_a+\bar{h}_{a}(g_a)$ is a critical point of $\bar{\Psi}_{a, K}$.

(iii) If $g_a\in \bar{G}$ and $\bar{H}_a$ is $C^k$ with $k\ge 2$ in a neighborhood of the trajectory of
$g_a+\bar{h}_{a}(g_a)$, then $\bar{\psi}_{a,K}$ is $C^{k-1}$ in a neighborhood of $g_a$.
In particular, if $g_a$ is a nonzero critical point
of $\bar{\psi}_{a,K}$, then $\bar{\psi}_{a,K}$ is $C^2$ in a neighborhood of the critical orbit $S^1\cdot g_a$. The
index and nullity of $\bar{\Psi}_{a,K}$ at $g_a+\bar{h}_{a}(g_a)$ defined in Definition 2.7 coincide
with the Morse index and nullity of $\bar{\psi}_{a,K}$ at $g_a$.

(iv) For any $\kappa\in\R$, we denote by¡¤¡¤
\bea  \wtd{\bar{\Lambda}}_{a,K}^\kappa=\{g\in \bar{G} \;|\; \bar{\psi}_{a,K}(g)\le\kappa\}.   \lb{2.39}\eea
Then the natural embedding $\wtd{\bar{\Lambda}}_{a, K}^\kappa \hookrightarrow {\bar{\Lambda}}_{a, K}^\kappa $ given by
$g\mapsto g+\bar{h}_{a}(g)$ is an $S^1$-equivariant homotopy equivalence.

(v) The functionals $a\mapsto\bar{\psi}_{a,K}$ is continuous in $a$ in the $C^1$ topology. Moreover
$a\mapsto\bar{\psi}^{\prime\prime}_{a,K}$ is continuous in a neighborhood of the critical orbit $S^1\cdot g_a$.}

{\bf Proof.} Let $x(t)=e^{JLt}x_0$ for some $\frac{LT}{2\pi}\in 2\Z-1$
and $x_0\in\R^{2n}$, then $-J\dot{x}+Kx=(L+K)x$. Thus $\{-\frac{1}{L+K}\mid \frac{LT}{2\pi}\in 2\Z-1\}$ is the set
of all the eigenvalues of $-\bar{M}_K$. By the convexity of $\bar{H}^*_{a,K}$, we have
\be (\bar{H}_{a,K}^{*\prime}(u)-\bar{H}_{a,K}^{*\prime}(v), u-v)\ge
            \omega|u-v|^2,\quad \forall\;a\in[a_1, a_2],\; u, v\in\R^{2n}, \lb{2.40}\ee
for some $\omega>0$. Hence we can use the proof of Proposition 3.9 of \cite{Vit2} to obtain the subspace $\bar{G}$ and
the map $\bar{h}_a$. In fact, Let $\bar{G}$ be the subspace of $L_s^{2}({\bf R}/(T {\bf Z}); {\bf R}^{2n})$
generated by the eigenvectors of $-\bar{M}_K$ whose eigenvalues are less than $-\frac{\omega}{2}$, i.e.,
$$ \bar{G} = span\{e^{JLt}x_0\mid -\frac{1}{L+K}<-\frac{\omega}{2}, \frac{LT}{2\pi}\in 2\Z-1, x_0\in\R^{2n}\}, $$
and $\bar{h}_a(g)$ is defined by the equation
\be \frac{\partial}{\partial h}\bar{\Psi}_{a, K}(g+\bar{h}_a(g))=0.  \lb{2.41}\ee
Then (i)-(iii) follow from Proposition 3.9 of \cite{Vit2}, and (iv) follows from Lemma 5.1 of \cite{Vit2},
the proof of (v) is the same as that of Lemma 2.10 (v) of \cite{LLW1}.\hfill\hb

{\bf Proposition 2.12.} {\it For all $b\ge a>\frac{\tau}{T}$, let $\bar{F}_{b,K}$ be the functional defined by
(\ref{2.21}), and $x_b$ be the critical point of $\bar{F}_{b,K}$ so that $x_b$ corresponds to a fixed symmetric closed
characteristic $(\tau,y)$ on $\Sigma$ for all $b\ge a$. Then the index $i_{-1}(A_b,1)$ and nullity $\nu_{-1}(A_b,1)$
are constants for all $b\ge a$. In particular, when $\bar{H}_b$ is $\aa$-homogenous for some $\aa\in (1,2)$ near
the image set of $x_b$, the index and nullity coincide with those defined for the Hamiltonian
$\bar{H}(x)=j(x)^\alpha$ for all $x\in\R^{2n}$. Especially $1\le \nu_{-1}(A_b,1)\le 2n-1$ always holds.}

{\bf Proof.} Denote by $R(t)$ the fundamental solution of the linearized system (\ref{2.33}) satisfying $R(0)=I_{2n}$.
Then by Lemma 1.6.11 of \cite{Eke2}, whose proof does not need the convexity of $\Sigma$, we have
\be  R(t)T_{y(0)}\Sigma\subset T_{y(\tau t/T)}\Sigma.   \lb{2.42}\ee
Then the completely same argument of Proposition 2.9 of \cite{LLo1} proves that $\nu_{-1}(A_a,1)$ is constant
for all $\bar{H}_a$ satisfying Proposition 2.4 with $a>\frac{\tau}{T}$ and $1\le \nu_{-1}(A_a,1)\le 2n-1$.

For any $b>a>\frac{\tau}{T}$, by (iii) of Lemma 2.2, we can construct a continuous family of $\bar{\Psi}_{c, K}$ with
$c\in[a, b]$ such that $\bar{H}_b$ is homogenous of degree $\aa=\aa_b$ near the image set of $x_b$. Now we can use
Lemma 2.11 (v) to obtain a continuous family of $\bar{\psi}_{c,K}$ such that $\bar{\psi}_{c,K}^{\prime\prime}(g_c)$
depends continuously on $c\in [a,b]$, where $g_c$ is the critical point of $\bar{\psi}_{c,K}$ corresponding to
$\bar{M}_K^{-1}{x}_c$. Because $\dim\ker\bar{\psi}_{c,K}^{\prime\prime}(g_c)=\bar{\nu}_K(\bar{M}_K^{-1}{x}_c)
=\nu_{-1}(A_c,1)={\rm constant}$, the index of $\bar{\psi}_{c,K}^{\prime\prime}(g_c)=\bar{i}_K(\bar{M}_K^{-1}{x}_c)
=i_{-1}(A_c,1)+e(K)$ must be constant too. Thus $i_{-1}(A_b,1)$ is
constant for all $b\ge a$. Note that here we used (\ref{2.29}), (\ref{2.35}), and Lemma 2.11 (iii).
Since the index $i_{-1}(A_b,1)$ and nullity $\nu_{-1}(A_b,1)$ only depend on the
value of $\bar{H}_b$ near the image set of $x_b$, then the index and nullity
coincide with those defined for the Hamiltonian $\bar{H}(x)=j(x)^\alpha$, $\forall~x\in\R^{2n}$.
The proof is complete.\hfill\hb

By Proposition 4.1 of \cite{Vit2} and Lemma 2.11, using the same proof of Proposition 2.12 of \cite{LLW1}, we also have:

{\bf Proposition 2.13.} {\it $\bar{\Psi}_{a, K}$ satisfies the Palais-Smale condition on $\bar{E}$,
and $\bar{F}_{a, K}$ satisfies
the Palais-Smale condition on $\bar{X}$, when $\frac{KT}{2\pi}\notin 2{\bf Z}-1$.}

Now for a critical point $u$ of $\bar{\Psi}_{a, K}$ and the corresponding $x=\bar{M}_K u$ of $\bar{F}_{a, K}$, let
\bea
\bar{\Lm}_{a,K}(u) &=& \bar{\Lm}_{a,K}^{\bar{\Psi}_{a, K}(u)}
   = \{w\in L^{2}_s(\R/(T\Z), \R^{2n}) \;|\; \bar{\Psi}_{a, K}(w)\le\bar{\Psi}_{a,K}(u)\},  \lb{2.43}\\
\bar{X}_{a,K}(x) &=& \bar{X}_{a,K}^{\bar{F}_{a,K}(x)} =
\{y\in W_s^{1, 2}(\R/(T\Z), \R^{2n}) \;|\; \bar{F}_{a,K}(y)\le \bar{F}_{a,K}(x)\}. \lb{2.44}\eea
Then both sets are $S^1$-invariant. Denote by $\crit(\bar{\Psi}_{a, K})$ the set of critical points of
$\bar{\Psi}_{a, K}$. Because $\bar{\Psi}_{a,K}$ is $S^1$-invariant,
$S^1\cdot u$ becomes a critical orbit if $u\in \crit(\bar{\Psi}_{a, K})$.
Note that by the condition $(F^\prime)$, Lemma 2.3, Proposition 2.4 and Lemma 2.5,
the number of critical orbits of $\bar{\Psi}_{a, K}$
is finite. Hence as usual we can make the following definition.

{\bf Definition 2.14.} {\it Suppose $u$ is a nonzero critical point of $\bar{\Psi}_{a, K}$, and $\Nn$ is
an $S^1$-invariant open neighborhood of $S^1\cdot u$ such that
$\crit(\bar{\Psi}_{a,K})\cap (\bar{\Lm}_{a,K}(u)\cap \Nn)
= S^1\cdot u$. Then the $S^1$-critical modules of $S^1\cdot u$ are defined by
\begin{eqnarray} C_{S^1,q}(\bar{\Psi}_{a}, S^1\cdot
u)=H_q((\bar{\Lambda}_{a, K}(u)\cap\mathcal {N})_{S^1},
((\bar{\Lambda}_{a, K}(u)\backslash S^1\cdot u)\cap\mathcal
{N})_{S^1}),\nn\end{eqnarray}
Similarly, we define the $S^1$-critical modules $C_{S^1,\; q}(\bar{F}_{a, K}, \;S^1\cdot x)$ of $S^1\cdot x$
for $\bar{F}_{a, K}$.}

By the same argument of Proposition 3.2 of \cite{LLW1}, we have the following for critical modules.

{\bf Proposition 2.15.} {\it  For any
$\frac{\tau}{T}<a_1<a_2<+\infty$, let $K$ be a fixed sufficiently large real number so that (\ref{2.19}) holds
for all $a\in [a_1, a_2]$. Then the critical module $C_{S^1,\; q}(\bar{F}_{a, K}, \;S^1\cdot x)$ is
independent of the choice of $\bar{H}_a$ defined in Proposition 2.4 for any $a\in [a_1, a_2]$ in the sense
that if $x_i$ is a solution of (\ref{2.17})
with Hamiltonian function $\bar{H}_{a_i}(x)$ with $i=1$ and $2$ respectively such that both $x_1$ and $x_2$
correspond to the same symmetric closed characteristic $(\tau,y)$ on $\Sigma$, then we have
\be C_{S^1,\;q}(\bar{F}_{a_1,K},\;S^1\cdot {x}_1) \cong C_{S^1,\;q}(\bar{F}_{a_2,K},\;S^1\cdot{x}_2), \qquad
\forall\,q\in \Z. \lb{2.45}\ee In other words, the critical modules are independent of the choices of all
$a>\frac{\tau}{T}$, the function $\vf_a$ satisfying
(i)-(ii) of Lemma 2.2, and $\bar{H}_a$ satisfying Proposition 2.4. }

Now we fix an $a>\frac{\hat{\tau}}{T}$, and write $\bar{F}_K$ and $\bar{H}$ for $\bar{F}_{a, K}$ and $\bar{H}_a$
respectively. We suppose also that $K\in {\bf R}$ satisfies (\ref{2.19}), i.e.,
\be  \bar{H}_K(x) = \bar{H}(x)+\frac{1}{2}K|x|^2 \qquad {\rm is\;strictly\;convex}. \lb{2.46}\ee

By Lemma 2.5, the critical points of $\bar{F}_{K}$ which are solutions of (\ref{2.17}) are the same for any
$K$ satisfying $\frac{KT}{2\pi}\notin 2{\bf Z}-1$. Recall that $e(K)=2n([K\frac{T}{2\pi}]-[K\frac{T}{4\pi}])$
in Theorem 2.8. The same argument of Theorem 3.3 of \cite{LLW1} proves

{\bf Theorem 2.16.} {\it Suppose $\bar{x}$ is a nonzero critical point of $\bar{F}_{K}$.
Then the $S^1$-critical module $C_{S^1,e(K)+l}(\bar{F}_{K}, S^1\cdot\bar{x})$ is independent
of the choice of $K$ for $\frac{KT}{2\pi}\notin 2{\bf Z}-1$, i.e.,
\be C_{S^1,e(K)+l}(\bar{F}_{K}, S^1\cdot\bar{x}) \cong C_{S^1,e(K^\prime)+l}(\bar{F}_{K^\prime}, S^1\cdot\bar{x}),  \lb{2.47}\ee
where $\frac{KT}{2\pi}$, $\frac{K^\prime T}{2\pi}\notin 2{\bf Z}-1$, $l\in \Z$, and
both $K$ and $K^\prime$ satisfy (\ref{2.46}).}

Now we fix $a$ and let $u_K\neq 0$ be a critical point of $\bar{\Psi}_{a, K}$ with multiplicity $\mul(u_K)=m$,
that is, $u_K$ corresponds to a symmetric closed characteristic $(\tau, y)\subset\Sigma$ with $(\tau, y)$
being $m$-iteration of some prime symmetric closed characteristic, where $m$ is odd.
Precisely, by Proposition 2.4 and Lemma 2.5, we have $u_K=-J\dot x+Kx$ with $x$
being a solution of (\ref{2.17}) and $x=\rho y(\frac{\tau t}{T})$ with $\frac{\vf_a^\prime(\rho)}{\rho}=\frac{\tau}{aT}$.
Moreover, $(\tau, y)$ is a symmetric closed characteristic on $\Sigma$ with minimal period $\frac{\tau}{m}$.
Hence the isotropy group satisfies $\{\theta\in S^1\;|\;\theta\cdot u_K=u_K\}=\Z_m$ and the orbit of $u_K$, namely,
$S^1\cdot u_K\cong S^1/\Z_m\cong S^1$. By Lemma 2.11, we obtain a critical point $g_K$ of $\bar{\psi}_{a,K}$
corresponding to $u_K$, and then its isotropy group satisfies $\{\theta\in S^1\;|\;\theta\cdot g_K=g_K\}=\Z_m$.
Let $p: N(S^1\cdot g_K)\rightarrow S^1\cdot g_K$ be the normal bundle of $S^1\cdot g_K$ in $\bar{G}$ (as defined
in Lemma 2.11) and let $p^{-1}(\theta\cdot g_K)=N(\theta\cdot g_K)$ be the fibre over $\theta\cdot g_K$,
where $\theta\in S^1$. Let $DN(S^1\cdot g_K)$ be the $\varrho$ disk bundle of $N(S^1\cdot g_K)$ for
some $\varrho>0$ sufficiently small, i.e.,
$DN(S^1\cdot g_K)=\{\xi\in N(S^1\cdot g_K)\;| \; \|\xi\|<\varrho\}$ which is identified by the exponential map with a
subset of $\bar{G}$, and let $DN(\theta\cdot g_K)=p^{-1}(\theta\cdot g_K)\cap DN(S^1\cdot g_K)$ be the disk over
$\theta\cdot g_K$. Clearly, $DN(\theta\cdot g_K)$ is $\Z_m$-invariant and we have
$DN(S^1\cdot g_K)=DN(g_K)\times_{\Z_m}S^1$, where the $Z_m$ action is given by
$$(\theta, v, t)\in \Z_m\times DN(g_K)\times S^1\mapsto (\theta\cdot v, \;\theta^{-1}t)\in DN(g_K)\times S^1. $$
Hence for an $S^1$ invariant subset $\Gamma$ of $DN(S^1\cdot g_K)$, we have
$\Gamma/S^1=(\Gamma_{g_K}\times_{\Z_m}S^1)/S^1=\Gamma_{g_K}/\Z_m$, where $\Gamma_{g_K}=\Gamma\cap DN(g_K)$.

For a $\Z_m$-space pair $(A, B)$, let
\bea H_{\ast}(A, B)^{\pm\Z_m} = \{\sigma\in H_{\ast}(A, B)\,|\,L_{\ast}\sigma=\pm \sigma\},\nn\eea
where $L$ is a generator of the $\Z_m$-action, we have

{\bf Lemma 2.17.} {\it Suppose $u_K\neq 0$ is a critical point of $\bar{\Psi}_{a, K}$ with $\mul(u_K)=m$,
$g_K$ is a critical
point of $\bar{\psi}_{a,K}$ corresponding to $u_K$. Then we have}
\bea C_{S^1,\;\ast}(\bar{\Psi}_{a,K}, \;S^1\cdot u_K)
\cong H_\ast((\wtd{\bar{\Lambda}}_{a,K}(g_K)\cap DN(g_K)),\;
    ((\wtd{\bar{\Lambda}}_{a,K}(g_K)\setminus\{g_K\})\cap DN(g_K)))^{\Z_m},   \nn\eea
where $\wtd{\bar{\Lambda}}_{a,K}(g_K)=\{g\in \bar{G} \;|\; \bar{\psi}_{a, K}(g)\le\bar{\psi}_{a, K}(g_K)\}$.

{\bf Proof.} We replace Lemma 2.10 used in the proof of Lemma 4.1 of \cite{LLW1} by the above Lemma 2.11, then
the proof follows from that of Lemma 4.1 of \cite{LLW1}.\hfill\hb

\setcounter{equation}{0}
\section{Periodic property of critical modules for symmetric closed characteristics}%Section 3

In this section, we use Lemmas 2.11 and 2.17 to obtain the periodic property of critical modules.

By (\ref{2.21}) and (\ref{2.24}), we have
$C_{S^1,\; q}(\bar{\Psi}_{a, K}, \;S^1\cdot u_K)\cong C_{S^1,\; q}(\bar{F}_{a, K}, \;S^1\cdot x)$.
By Proposition 2.15, the module $C_{S^1,\; q}(\bar{F}_{a, K}, \;S^1\cdot x)$ is independent of the choice of
the Hamiltonian function $\bar{H}_a$ whenever $\bar{H}_a$ satisfies conditions in Proposition 2.4.
Hence in order to compute the critical modules, we can
choose $\bar{\Psi}_{a, K}$ with $\bar{H}_a$ being positively homogeneous of degree $\aa=\aa_a$ near the
image set of every nonzero solution $x$ of (\ref{2.17}) corresponding to some symmetric closed characteristic
$(\tau,y)$ with period $\tau$ being strictly less than $aT$.

In other words, for a given $a>0$, we choose $\vth\in (0,1)$ first such that
$[aT\vth, aT(1-\vth)]\supset \per(\Sigma)\cap (0,aT)$ holds by the definition of the set $\per(\Sg)$
and the assumption (F). Then we choose $\aa=\aa_a\in (1,¡¢£¬2)$ sufficiently close to $2$ by (iii) of Lemma 2.2
such that $\vf_a(t)=ct^\aa$ for some constant $c>0$ and $\aa\in(1,\,2)$ whenever
$\frac{\vf_a^\prime(t)}{t}\in [\vth, 1-\vth]$. Now we suppose that
$\vf_a$ satisfies (iii) of Lemma 2.2.

Now we consider iterations of critical points of $\bar{\Psi}_{a, K}$. Suppose $u_K\neq 0$ is a critical point of
$\bar{\Psi}_{a,K}$ with $\mul(u_K)=m$, where $m$ is odd, and $g_K$ is the critical point of $\bar{\psi}_{a,K}$
corresponding to $u_K$. By Proposition 2.4 and Lemma 2.5, we have $u_K=-J\dot x+Kx$ with $x$ being a solution
of (\ref{2.17}) and $x=\rho y(\frac{\tau t}{T})$ with $\frac{\vf_a^\prime(\rho)}{\rho}=\frac{\tau}{aT}$.
Moreover, $(\tau, y)$ is a symmetric closed characteristic on $\Sigma$ with minimal period $\frac{\tau}{m}$.
For any $p\in 2\N-1$ satisfying $p\tau<aT$,
we choose $K$ such that $\frac{pKT}{2\pi}\notin 2{\bf Z}-1$, then the $p$th iteration $u_{pK}^p$ of $u_K$ is
given by $-J\dot x^p+pKx^p$, where $x^p$ is the unique solution of (2.17) corresponding to $(p\tau, y)$ and
is a critical point of $\bar{F}_{a, pK}$, that
is, $u_{pK}^p$ is the critical point of $\bar{\Psi}_{a, pK}$ corresponding to $x^p$. Hence we have
\bea
&& x(t)=\left(\frac{\tau}{c\alpha a}\right)^\frac{1}{\alpha-2}y(\tau t), \quad
      x^p(t)=\left(\frac{p\tau}{c\alpha a}\right)^\frac{1}{\alpha-2}y(p\tau t)=p^{\frac{1}{\alpha-2}}x(pt), \nn\\
&& u_K(t)=-J\dot {x}(t)+Kx(t), \quad u_{pK}^p(t)=-J\dot x^p(t)+pKx^p(t)=p^{\frac{\alpha-1}{\alpha-2}}u_K(pt). \nn\eea

We define the $p$th iteration $\phi^p$ on $L_s^{2}(\R/(T\Z); {\bf R}^{2n})$ by
\be \phi^p: v_K(t)\mapsto v^p_{pK}(t)\equiv p^\frac{\alpha-1}{\alpha-2}v_K(pt). \lb{3.1}\ee

Now we use the notations in Lemma 2.11, we choose $\bar{G}_{pK}$ in Lemma 2.11 for $\bar{\Psi}_{a, pK}$ such that
$\bar{G}_{pK}\supseteq \phi^p(\bar{G}_K)$, where we write $\bar{G}_K$ for $\bar{G}$ to indicate its dependence on $K$.
For $g\in DN(S^1\cdot g_K)$, as (4.12) of \cite{LLW1}, we have
\bea \bar{\psi}_{a, pK}(\phi^p(g))=p^{\frac{\alpha}{\alpha-2}}\bar{\psi}_{a, K}(g).  \lb{3.2}\eea

We define a new inner product $\<\cdot,\cdot\>_p$ on $L_s^2(\R/\Z, \R^{2n})$ by
\be \<v,w\>_p = p^\frac{2(\alpha-1)}{2-\alpha}\<v, w\>.   \lb{3.3}\ee
Then $\phi^p: DN(g_K)\rightarrow DN(g_{pK}^p)$ is an isometry from the standard inner product to the above one,
where $g_{pK}^p=\phi^p(g_K)$ is the critical point of $\bar{\psi}_{a,pK}$ corresponding to $u_{pK}^p$ and the radii of
the two normal disk bundles are suitably chosen. Clearly $\phi^p(DN(g_K))$ consists of points in $DN(g_{pK}^p)$
which are fixed by the $\Z_p$-action. Since the $\Z_p$-action on $DN(g_{pK}^p)$ is an isometry and
$f\equiv\bar{\psi}_{a, pK}|_{DN(g_{pK}^p)}$ is $\Z_p$-invariant, we have
\be f^{\prime\prime}(g) =\left(\matrix{(f|_{\phi^p(DN(g_K))})^{\prime\prime} \quad 0\cr \qquad
             0\qquad\qquad\;\;\ast}\right),\quad \forall g\in \phi^p(DN(g_K)).  \lb{3.4}\ee
Moreover, we have
\be f^\prime(g)=(f|_{\phi^p(DN(g_K))})^{\prime},\quad \forall g\in \phi^p(D N(g_K)).  \lb{3.5}\ee
Now we can apply the results by D. Gromoll and W. Meyer \cite{GrM1} to the manifold $DN(g_{pK}^p)$ with $g_{pK}^p$
as its unique critical point. Then $\mul(g_{pK}^p)=pm$ is the multiplicity of $g_{pK}^p$ and the isotropy group
$\Z_{pm}\subseteq S^1$ of $g_{pK}^p$ acts on $DN(g_{pK}^p)$ by isometries. According to Lemma 1 of \cite{GrM1}, we
have a $\Z_{pm}$-invariant decomposition of $T_{g_{pK}^p}(DN(g_{pK}^p))$
\be  T_{g_{pK}^p}(DN(g_{pK}^p)) = \bar{V}^+\oplus \bar{V}^-\oplus \bar{V}^0 = \{(x_+, x_-, x_0)\}  \lb{3.6}\ee
with $\dim \bar{V}^-=\bar{i}(g_{pK}^p)=\bar{i}_{pK}(u_{pK}^p)$,
$\dim \bar{V}^0=\bar{\nu}(g_{pK}^p)-1=\bar{\nu}_{pK}(u_{pK}^p)-1$ (cf. Lemma 2.11(iii)),
and a $\Z_{pm}$-invariant neighborhood $\bar{B}=\bar{B}_+\times \bar{B}_-\times \bar{B}_0$ for $0$
in $T_{g_{pK}^p}(DN(g_{pK}^p))$ together
with two $\Z_{pm}$-invariant diffeomorphisms
$$ \bar{\Phi} : \bar{B}=\bar{B}_+\times \bar{B}_-\times \bar{B}_0\to \bar{\Phi}(\bar{B}_+
\times \bar{B}_-\times \bar{B}_0)\subset DN(g_{pK}^p),  $$
and
$$ \bar{\eta} : \bar{B}_0\to W(g_{pK}^p)\equiv\bar{\eta}(\bar{B}_0)\subset DN(g_{pK}^p), $$
and $\bar{\Phi}(0)=\bar{\eta}(0)=g_{pK}^p$, such that
\be \bar{\psi}_{a,pK}\circ\bar{\Phi}(x_+,x_-,x_0)=|x_+|^2 - |x_-|^2 +
\bar{\psi}_{a,pK}\circ\bar{\eta}(x_0), \lb{3.7}\ee
with $d(\bar{\psi}_{a, pK}\circ \bar{\eta})(0)=d^2(\bar{\psi}_{a, pK}\circ\bar{\eta})(0)=0$.
As usual, we call $W(g_{pK}^p)$ a local
characteristic manifold, and $U(g_{pK}^p)=\bar{B}_-$ a local negative disk at $g_{pK}^p$. By the proof of Lemma 1 of
\cite{GrM1}, $W(g_{pK}^p)$ and $U(g_{pK}^p)$ are $\Z_{pm}$-invariant. It follows from (\ref{3.7}) that $g_{pK}^p$
is an isolated critical point of $\bar{\psi}_{a, pK}|_{DN(g_{pK}^p)}$. Then as in Lemma 6.4 of \cite{Rad1}, we have
\bea
&& H_\ast(\wtd{\bar{\Lambda}}_{a, pK}(g_{pK}^p)\cap DN(g_{pK}^p),\;
     (\wtd{\bar{\Lambda}}_{a, pK}(g_{pK}^p)\setminus\{g_{pK}^p\})\cap DN(g_{pK}^p))\nn\\
&&\quad = \bigoplus_{q\in\Z}H_q (U(g_{pK}^p),U(g_{pK}^p)\setminus\{g_{pK}^p\}) \nn\\
&&\qquad\qquad\qquad \otimes H_{\ast-q}(W(g_{pK}^p)\cap \wtd{\bar{\Lambda}}_{a,pK}(g_{pK}^p),
       (W(g_{pK}^p)\setminus\{g_{pK}^p\})\cap \wtd{\bar{\Lambda}}_{a,pK}(g_{pK}^p)), \lb{3.8}\eea
where
\be  H_q(U(g_{pK}^p),U(g_{pK}^p)\setminus\{g_{pK}^p\} )
    = \left\{\matrix{\Q, & {\rm if\;}q=\bar{i}_{pK}(u_{pK}^p),  \cr
                      0, & {\rm otherwise}. \cr}\right.    \lb{3.9}\ee

Now we have the following proposition.

{\bf Proposition 3.1.} {\it For any $p\in2\N-1$, we choose $K$ such that $\frac{pKT}{2\pi}\notin 2{\bf Z}-1$.
Let $u_K\neq 0$ be a critical point of $\bar{\Psi}_{a, K}$ with $\mul(u_K)=1$, $u_K=-J\dot x+Kx$ with $x$ being
a critical point of $\bar{F}_{a, K}$. Then for all $q\in\Z$, we have
\bea
&& C_{S^1,\; q}(\bar{\Psi}_{a,pK},\;S^1\cdot u_{pK}^p) \nn\\
&&\quad\cong \left(\frac{}{}H_{q-\bar{i}_{pK}(u_{pK}^p)}(W(g_{pK}^p)\cap \wtd{\bar{\Lambda}}_{a,pK}
(g_{pK}^p),(W(g_{pK}^p)\setminus\{g_{pK}^p\})\cap \wtd{\bar{\Lambda}}_{a,pK}(g_{pK}^p))\right)^{\Z_p}. \lb{3.10}\eea
In particular, if $u_{pK}^p$ is
non-degenerate, i.e., $\bar{\nu}_{pK}(u_{pK}^p)=1$, then}
\be C_{S^1,\; q}(\bar{\Psi}_{a,pK},\;S^1\cdot u_{pK}^p)
    = \left\{\matrix{\Q, & {\rm if\;}q=\bar{i}_{pK}(u_{pK}^p),  \cr
                      0, & {\rm otherwise}. \cr}\right. \lb{3.11}\ee

{\bf Proof.} Suppose $\theta$ is a generator of the linearized $\Z_p$-action on $U(g_{pK}^p)$. Then $\theta(\xi)=\xi$
if and only if $\xi\in T_{g_{pK}^p}(\phi^p(DN(g_K)))$. Hence it follows from (\ref{3.2}) and (\ref{3.4}) that
$\xi=(\phi^p)_\ast(\xi^\prime)$ for a unique $\xi^\prime\in T_{g_K}(DN(g_K))^-$. Hence the proof of Satz 6.11 in
\cite{Rad1}, Proposition 2.8 in \cite{BaL1} yield this proposition. Note that
$\bar{i}_{pK}(u_{pK}^p)=2n([pKT/{2\pi}]-[pKT/{4\pi}])+i_{-1}(A_a, p)$ and
$\bar{i}_{K}(u_{K})=2n([KT/{2\pi}]-[KT/{4\pi}])+i_{-1}(A_a, 1)$ follow from
Theorem 2.8, and when $p$ is odd, $i_{-1}(A_a, p)-i_{-1}(A_a, 1)$ is always even (cf.
Theorem 9.2.1, Theorem 9.3.4 of \cite{Lon3}).\hfill\hb

{\bf Definition 3.2.} {\it For any $p\in2\N-1$, we choose $K$ such that $\frac{pKT}{2\pi}\notin 2{\bf Z}-1$.
Let $u_K\neq 0$ be
a critical point of $\bar{\Psi}_{a,K}$ with $\mul(u_K)=1$, $u_K=-J\dot x+Kx$ with $x$ being a critical point
of $\bar{F}_{a, K}$.
Then for all $l\in\Z$, let
\bea
\bar{k}_{l}(u_{pK}^p) &=& \dim\left(\frac{}{}H_{l}(W(g_{pK}^p)\cap
  \wtd{\bar{\Lambda}}_{a,pK}(g_{pK}^p),(W(g_{pK}^p)
  \bs\{g_{pK}^p\})\cap \wtd{\bar{\Lambda}}_{a,pK}(g_{pK}^p))\right)^{\Z_p}.
          \lb{3.12} \eea
$\bar{k}_l(u_{pK}^p)$'s are called critical type numbers of $u_{pK}^p$. }

{\bf Remark 3.3.} (i) Since
\bea C_{S^1,\;l+\bar{i}_{pK}(u_{pK}^p)}(\bar{\Psi}_{a,pK}, \;S^1\cdot u_{pK}^p)
&\cong& C_{S^1,\; l+\bar{i}_{pK}(x^p)}(\bar{F}_{a,pK}, \;S^1\cdot x^p)  \nn\\
&\cong& C_{S^1,\; l+e(pK)+i_{-1}(A_a, p)}(\bar{F}_{a,pK}, \;S^1\cdot x^p),  \nn\eea
by Theorem 2.16, we obtain that $\bar{k}_l(u_{pK}^p)$ is independent of the choice of $K$
and denote it by $\bar{k}_l(x^p)$,
here $\bar{k}_l(x^p)$'s are called critical type numbers of $x^p$.

(ii) By Proposition 2.12, we have $\bar{k}_{l}(u_{pK}^p)=0$ if $l\notin [0, 2n-2]$.

{\bf Lemma 3.4.} {\it Let $u_K\neq 0$ be a critical point of $\bar{\Psi}_{a, K}$ with $\mul(u_K)=1$. Suppose
$\bar{\nu}_{mK}(u_{mK}^m)=\bar{\nu}_{pmK}(u_{pmK}^{pm})$ for some $m, p\in2\N-1$. Then we have
$\bar{k}_{l}(u_{mK}^m)=\bar{k}_{l}(u_{pmK}^{pm})$ for all $l\in\Z$.}

{\bf Proof.} From the above arguments, we obtain our lemma by the same method of the proof of
Lemma 4.5 of \cite{LLW1}.\hfill\hb

{\bf Proposition 3.5.} {\it Let $x\neq 0$ be a critical point of $\bar{F}_{a,K}$ with $\mul(x)=1$ corresponding to a
critical point $u_K$ of $\bar{\Psi}_{a, K}$. Then there exists a minimal $\bar{K}(x)\in 2\N$ such that
\bea
&& \nu_{-1}(x^{p+\bar{K}(x)})=\nu_{-1}(x^p),\quad i_{-1}(x^{p+\bar{K}(x)})-i_{-1}(x^p)\in 2\Z,  \lb{3.13}\\
&& \bar{k}_l(x^{p+\bar{K}(x)})=\bar{k}_l(x^p), \qquad\forall p\in 2\N-1,\;l\in\Z, \lb{3.14}\eea
where $i_{-1}(x^p)$ is defined in Theorem 2.8.
We call $\bar{K}(x)$ the minimal period of critical modules of iterations of the functional $\bar{F}_{a, K}$ at $x$.}

{\bf Proof.} We replace Lemma 3.3 used in the proof of Proposition 3.4 of \cite{LLo1} by the above Lemma 3.4, then
our proposition follows by Proposition 3.4 of \cite{LLo1}.\hfill\hb

In the following, we give the definitions of indices and Euler characteristics for symmetric closed characteristics.
Let $\bar{F}_{a,K}$ be any function defined by (\ref{2.21}) with $\bar{H}_a$ satisfying
Proposition 2.4, we do not require $\bar{H}_a$ to be homogeneous near its critical points.

{\bf Definition 3.6.} {\it Suppose the condition $(F^\prime)$ in Section 2 holds. For every symmetric closed
characteristic $(\tau,y)$ on $\Sigma$, let $aT>\tau$ and choose $\vf_a$ to satisfy (i)-(ii) of Lemma 2.2.
Determine $\rho$ uniquely by $\frac{\vf_a'(\rho)}{\rho}=\frac{\tau}{aT}$. Let $x=\rho y(\frac{\tau t}{T})$.
Then we define the index $\bar{i}(\tau,y)$ and nullity $\bar{\nu}(\tau,y)$ of $(\tau,y)$ by (cf. Theorem 2.8)
$$ \bar{i}(\tau,y)=i_{-1}(x), \qquad \bar{\nu}(\tau,y)=\nu_{-1}(x). $$
Then the mean index of $(\tau, y)$ is defined by }
\be \hat{\bar{i}}(\tau,y) = \lim_{m\rightarrow\infty}\frac{\bar{i}((2m-1)\tau, y)}{2m-1}.  \lb{3.15}\ee

Note that by Proposition 2.12, the index and nullity are well defined and are independent of the
choice of $aT>\tau$ and $\vf_a$ satisfying (i)-(ii) of Lemma 2.2.

For a prime symmetric closed characteristic $(\tau, y)$ on $\Sigma$, we denote simply by
$y^m\equiv(m\tau, y)$ for $m\in2\N-1$. By Proposition 2.15, we can define the critical type
numbers $\bar{k}_l(y^m)$ of $y^m$ to be $\bar{k}_l(x^m)$, where $x^m$
is the critical point of $\bar{F}_{a, K}$ corresponding to $y^m$. We also define $\bar{K}(y)=\bar{K}(x)$,
where $\bar{K}(x)\in 2\N$ is given by Proposition 3.5.
Suppose $\Nn$ is an $S^1$-invariant open neighborhood of $S^1\cdot x^m$ such that
$\crit(\bar{F}_{a, K})\cap(\bar{X}_{a,K}(x^m)\cap\Nn)=S^1\cdot x^m$. Then we make the following definition

{\bf Definition 3.7.} {\it The Euler characteristic $\bar{\chi}(y^m)$ of $y^m$ is defined by
\bea \bar{\chi}(y^m)
&\equiv& \chi((\bar{X}_{a,K}(x^m)\cap\Nn)_{S^1},\; ((\bar{X}_{a,K}(x^m)\setminus S^1\cdot x^m)\cap\Nn)_{S^1}) \nn\\
&\equiv& \sum_{q=0}^{\infty}(-1)^q\dim C_{S^1,\;q}(\bar{F}_{a,K},\;S^1\cdot x^m). \lb{3.16}\eea
Here $\chi(A,B)$ denotes the usual Euler characteristic of the space pair $(A,B)$.
The average Euler characteristic $\hat{\bar{\chi}}(y)$ of $y$ is defined by}
\be \hat{\bar{\chi}}(y)=\lim_{N\to\infty}\frac{1}{N}\sum_{1\le m\le N}\bar{\chi}(y^{2m-1}). \lb{3.17}\ee

Note that by Proposition 2.15 and Theorem 2.16, $\bar{\chi}(y^m)$ is well defined and is independent of the
choice of $a$ and $K$. In fact, by Remark 3.3 (i), we have
\be  \bar{\chi}(y^m)=\sum_{l=0}^{2n-2}(-1)^{\bar{i}(y^m)+l}\bar{k}_l(y^m).  \lb{3.18}\ee
The following remark shows that $\hat{\bar{\chi}}(y)$ is well-defined and is a rational number.

{\bf Remark 3.8.} By (\ref{3.13}), (\ref{3.18}) and Proposition 3.5, we have
\bea \hat{\bar{\chi}}(y)
&=&\lim_{N\rightarrow\infty}\frac{1}{N}  \sum_{1\le m\le N\atop 0\le l\le 2n-2}(-1)^{\bar{i}(y^{2m-1})+l}
\bar{k}_l(y^{2m-1}) \nn\\
&=&\lim_{s\rightarrow\infty}\frac{2}{s\bar{K}(y)}\sum_{1\le m\le \bar{K}(y)/2,\; 0\le l\le 2n-2\atop 0\le p< s}
              (-1)^{\bar{i}(y^{p\bar{K}(y)+2m-1})+l}\bar{k}_l(y^{p\bar{K}(y)+2m-1}) \nn\\
&=&\frac{2}{\bar{K}(y)}\sum_{1\le m\le \bar{K}(y)/2\atop 0\le l\le 2n-2}(-1)^{\bar{i}(y^{2m-1})+l}\bar{k}_l(y^{2m-1}).
\lb{3.19}\eea
Therefore $\hat{\bar{\chi}}(y)$ is well defined and is a rational number. In particular, if all $y^m$s are
non-degenerate, then $\bar{\nu}(y^m)=1$ for all $m\in2\N-1$. Hence the proof of Proposition 3.5 yields $\bar{K}(y)=2$.
By (\ref{3.11}), we have
$$ \bar{k}_l(y^m)
    = \left\{\matrix{1, & {\rm if\;\;}  l=0  \cr
                     0, & {\rm otherwise}. \cr}\right.  $$
Hence (\ref{3.19}) implies
\be \hat{\bar{\chi}}(y) = (-1)^{\bar{i}(y)}. \lb{3.20}\ee

\setcounter{equation}{0}
\section{Proof of Theorem 1.1}%Section 4

In this section, we give the proof of Theorem 1.1. Firstly, we consider the contribution
of the origin to the Morse series of the functional $\bar{F}_{a,K}$ on $W_s^{1,2}(\R/\Z;\R^{2n})$.
The same argument as in Theorem 5.1 of \cite{LLW1} and Theorem 7.1 of
\cite{Vit2} yields the following:

{\bf Theorem 4.1.} {\it Fix an $a>0$ such that $\per(\Sg)\cap (0,aT)\not=\emptyset$. Then there exists an
$\epsilon_0>0$ small enough such that for any $\epsilon\in (0,\epsilon_0]$ we have
\bea H_{S^1,\;q+e(K)}(\bar{X}_{a, K}^{\epsilon}, \;\bar{X}_{a, K}^{-\epsilon}) = 0,  \quad \forall q\in \mathring {I}, \nn\eea
if $I$ is an interval of $\Z$ such that $I\cap [\bar{i}(\tau, y), \bar{i}(\tau, y)+\bar{\nu}(\tau, y)-1]=\emptyset$
for all symmetric closed characteristics $(\tau,\, y)$ on $\Sigma$ with $\tau\ge aT$.}

Now we prove Theorem 1.1.

Let $\bar{F}_{a, K}$ be a functional defined by (\ref{2.21}) for some $a, K\in\R$ sufficiently large and
by Proposition 2.6, let $\epsilon>0$ be
small enough such that $[-\epsilon, 0)$ contains no critical values of $\bar{F}_{a, K}$. We consider the exact sequence
of the triple $(\bar{X}, \bar{X}^{-\epsilon}, X^{-b})$ (for $b$ large enough)
\bea \rightarrow H_{S^1, *}(\bar{X}^{-\ep}, \bar{X}^{-b})
&\to& H_{S^1, *}(\bar{X},\bar{X}^{-b})  \nn\\
&\to& H_{S^1, *}( \bar{X}, \bar{X}^{-\ep})\to H_{S^1, *-1}( \bar{X}^{-\ep}, \bar{X}^{-b})\to \cdots,  \lb{4.1}\eea
where $\bar{X}=W_s^{1, 2}(\R/{T\Z}; \R^{2n})$. The normalized Morse series of $\bar{F}_{a, K}$ in
$ \bar{X}^{-\ep}\setminus \bar{X}^{-b}$ is defined, as usual, by
\be  \bar{M}_a(t)=\sum_{q\ge 0,\;1\le j\le p} \dim C_{S^1,\;q}(\bar{F}_{a, K}, \;S^1\cdot \bar{v}_j)t^{q-e(K)},
\lb{4.2}\ee
where we denote by $\{S^1\cdot \bar{v}_1, \ldots, S^1\cdot \bar{v}_p\}$ the critical orbits of $\bar{F}_{a, K}$
with critical values less than $-\epsilon$. We denote by $t^{e(K)}\bar{H}_a(t)$ the Poincar\'e series of
$H_{S^1, *}(\bar{X}^{-\ep}, \bar{X}^{-b})$,
$\bar{H}_a(t)$ is a Laurent series, and we have the equivariant Morse inequality
\be  \bar{M}_a(t)-\bar{H}_a(t)=(1+t)\bar{R}_a(t),  \lb{4.3}\ee
where $\bar{R}_a(t)$ is a Laurent series with nonnegative coefficients.

On the other hand, by Corollary 2.10 the Poincar\'e series of $H_{S^1, *}(\bar{X}, \bar{X}^{-b})$ is
$t^{e(K)}(1/(1-t^2))$. The Poincar\'e series of $H_{S^1, *}(\bar{X}, \bar{X}^{-\ep})$ is $t^{e(K)}\bar{Q}_a(t)$,
according to Theorem 4.1, if we set $\bar{Q}_a(t)=\sum_{k\in \Z}{\bar{q}_kt^k}$, then
\be   \bar{q}_k=0 \qquad\qquad \forall\;k\in \mathring {I},  \lb{4.4}\ee
where $I$ is defined in Theorem 4.1. Now using (\ref{4.1}) and Proposition 1 in Appendix 2 of \cite{Vit2}, these
results yield
\be  \bar{H}_a(t)-\frac{1}{1-t^2}+\bar{Q}_a(t) = (1+t)\bar{S}_a(t),   \lb{4.5}\ee
with $\bar{S}_a(t)$ a Laurent series with nonnegative coefficients. Adding up (\ref{4.3}) and (\ref{4.5}) yields
\be  \bar{M}_a(t)-\frac{1}{1-t^2}+\bar{Q}_a(t) = (1+t)\bar{U}_a(t),   \lb{4.6}\ee
where $\bar{U}_a(t)=\sum_{i\in \Z}{\bar{u}_it^i}$ also has nonnegative coefficients.

Now truncate (\ref{4.6}) at the degrees $2C$ and $2N$, where we set $C$ equal to $2n^2$, and $2N>2C$, and
write $\bar{M}_a^{2N}(2C; t)$, $\bar{Q}_a^{2N}(2C; t)\cdots$ for the truncated series. Then from (\ref{4.6}) we infer
\bea
&& \bar{M}_a^{2N}(2C; t)-\sum_{h=C}^{N}{t^{2h}}+\bar{Q}_a^{2N}(2C;t)   \nn\\
&& \qquad = (1+t)\bar{U}_a^{2N-1}(2C; t)+t^{2N}\bar{u}_{2N}+t^{2C}\bar{u}_{2C-1}.  \lb{4.7}\eea
By (\ref{4.4}), and the fact that for $a$ large enough $\mathring {I}$ contains $[2C, 2N]$, indeed let $\alpha>0$
such that any prime symmetric closed characteristic $(\tau, y)$ with $\hat{\bar{i}}(y)\neq 0$ has
$|\hat{\bar{i}}(y)|>\alpha$, then if
$k\geq aT/\min{\{\tau_i\}}$, we have
$|\bar{i}(y^k)| \sim k|\hat{\bar{i}}(y)|\geq k\alpha\geq {a\alpha T}/\min{\{\tau_i\}}$,
which tends to infinity as $a\to\infty$. So $\bar{Q}_a^{2N}(2C; t)=0$, and (\ref{4.7}) can be written
\be  \bar{M}_a^{2N}(2C;t)-\sum_{h=C}^{N}{t^{2h}} =
(1+t)\bar{U}_a^{2N-1}(2C;t)+t^{2N}\bar{u}_{2N}+t^{2C}\bar{u}_{2C-1}.  \lb{4.8}\ee
Changing $C$ into $-C$, $N$ into $-N$, and counting terms with $-2N\leq i\leq -2C$, we obtain
\be  \bar{M}_a^{-2C}(-2N; t) = (1+t)\bar{U}_a^{-2C-1}(-2N;t)+t^{-2N}\bar{u}_{-2N-1}+t^{-2C}\bar{u}_{-2C}.  \lb{4.9}\ee
Denote by $\{x_1, \ldots, x_k\}$ the critical points of $\bar{F}_{a, K}$ corresponding to $\{y_1,\ldots, y_k\}$. Note
that $\bar{v}_1,\ldots,\bar{v}_p$ in (\ref{4.2}) are odd iterations of $x_1,\ldots,x_k$. Since
$C_{S^1,\;q}(\bar{F}_{a, K}, \;S^1\cdot x_j^m)$ can be non-zero only for
$q=e(K)+\bar{i}(y_j^m)+l$ with $0\le l\le 2n-2$,
by Propositions 2.12, 3.1 and Remark 3.3, the normalized Morse series (\ref{4.2}) becomes
\be  \bar{M}_a(t) = \sum_{1\le j\le k,\; 0\le l\le 2n-2 \atop 1\le 2m_j-1<aT/\tau_j}
\bar{k}_l(y_j^{2m_j-1})t^{\bar{i}(y_j^{2m_j-1})+l}
  \;\;= \sum_{1\le j\le k,\; 0\le l\le 2n-2 \atop 1 \le m_j\le \bar{K}_j/2,\; s\bar{K}_j+2m_j-1<aT/\tau_j}
             \bar{k}_l(y_j^{2m_j-1})t^{\bar{i}(y_j^{s\bar{K}_j+2m_j-1})+l}, \lb{4.10}\ee
where $\bar{K}_j=\bar{K}(y_j)$ and $s\in\N_0$. The last equality follows from Proposition 3.5.

Write $\bar{M}(t)=\sum_{h\in \Z}w_ht^h$, where $\bar{M}(t)$ denotes $\bar{M}_a(t)$ as $a$ tends to infinity.
Then we have
\be w_h\ = \sum_{1\le j\le k,\; 0\le l\le 2n-2 \atop 1 \le m\le \bar{K}_j/2}
              \bar{k}_l(y_j^{2m-1})\,^\#\{s\in\N_0\,|\,\bar{i}(y_j^{s\bar{K}_j+2m-1})+l=h\},
              \quad \forall\;2C\leq|h|\le 2N. \lb{4.11}\ee
Note that the right hand side of (\ref{4.10}) contains only those terms satisfying $s\bar{K}_j+2m_j-1<\frac{aT}{\tau_j}$.
Thus (\ref{4.11}) holds for $2C\leq |h|\leq 2N$ by (\ref{4.10}).

{\bf Claim 1.} {\it $w_h\le C_1$ for $2C\leq |h|\leq 2N$ with $C_1$ being independent of $a, K$}.

In fact, we have
\bea
^\#\{s\in\N_0 &|& \bar{i}(y_j^{s\bar{K}_j+2m-1})+l=h \}  \nn\\
&=& \;^\#\{s\in\N_0 \;| \; \bar{i}(y_j^{s\bar{K}_j+2m-1})+l=h,\;|\bar{i}(y_j^{s\bar{K}_j+2m-1})-
(s\bar{K}_j+2m-1)\hat{\bar{i}}(y_j)|\le 2n\} \nn\\
&\le& \;^\#\{s\in\N_0 \;| \;|h-l-(s\bar{K}_j+2m-1)\hat{\bar{i}}(y_j)|\le 2n\}  \nn\\
&=& \;^\#\left\{s\in\N_0 \; \left|\;\frac{}{}\right.\;h-l-2n-(2m-1)\hat{\bar{i}}(y_j)\le s\bar{K}_j\hat{\bar{i}}(y_j)
       \le h-l+2n-(2m-1)\hat{\bar{i}}(y_j)\right\}  \nn\\
&\le&\; \frac{4n}{\bar{K}_j|\hat{\bar{i}}(y_j)|}+2,  \lb{4.12}\eea
where the first equality follows from the fact that
\be  |\bar{i}(y_j^m)-m\hat{\bar{i}}(y_j)|\le 2n,\quad \forall m\in2\N-1,\; 1\le j\le k,  \lb{4.13}\ee
which follows from Theorem 9.2.1 and Theorems 10.1.2 of \cite{Lon3}, Note that
$\bar{i}(y_j^{s\bar{K}_j+2m-1})+l=h\in[2C, 2N]$ holds only when $\hat{\bar{i}}(y_j)>0$ and
$\bar{i}(y_j^{s\bar{K}_j+2m-1})+l=h\in[-2N, -2C]$ holds
only when $\hat{\bar{i}}(y_j)< 0$. Hence Claim 1 holds.

Next we estimate $\bar{M}_a^{2N}(2C; -1)$ and $\bar{M}_a^{-2C}(-2N; -1)$. By (\ref{4.11}) we obtain
\bea
&& \bar{M}_a^{2N}(2C; -1) = \sum_{h=2C}^{2N} w_h(-1)^h   \nn\\
&& = \sum_{1\le j\le k,\; 0\le l\le 2n-2 \atop 1 \le m\le \bar{K}_j/2}(-1)^{\bar{i}(y_j^{2m-1})+l}\bar{k}_l(y_j^{2m-1})
              \,^\#\{s\in\N_0 \,|\, 2C\le \bar{i}(y_j^{s\bar{K}_j+2m-1})+l\le 2N\}.  \lb{4.14}\eea
Here the second equality holds by (3.13). Similarly, we have
\bea
&& \bar{M}_a^{-2C}(-2N; -1) = \sum_{h=-2N}^{-2C} w_h(-1)^h   \nn\\
&&= \sum_{1\le j\le k,\; 0\le l\le 2n-2 \atop 1 \le m\le \bar{K}_j/2}(-1)^{\bar{i}(y_j^{2m-1})+l}\bar{k}_l(y_j^{2m-1})
              \,^\#\{s\in\N_0 \,|\, -2N\le \bar{i}(y_j^{s\bar{K}_j+2m-1})+l\le -2C\}.  \qquad \lb{4.15}\eea

{\bf Claim 2.} {\it There is a real constant $C_2>0$ independent of $a,K$ such that
\bea
\left|\bar{M}_a^{2N}(2C; -1)-\sum_{1\le j\le k,\;0\le l\le 2n-2 \atop 1 \le m\le \bar{K}_j/2, \hat{\bar{i}}(y_j)>0}
            (-1)^{\bar{i}(y_j^{2m-1})+l}\bar{k}_l(y_j^{2m-1})\frac{2N}{\bar{K}_j\hat{\bar{i}}(y_j)}\right|
               &\le& C_2,  \lb{4.16}\\
\left|\bar{M}_a^{-2C}(-2N; -1)-\sum_{1\le j\le k,\; 0\le l\le 2n-2 \atop 1 \le m\le \bar{K}_j/2, \hat{\bar{i}}(y_j)<0}
            (-1)^{\bar{i}(y_j^{2m-1})+l}\bar{k}_l(y_j^{2m-1})\frac{2N}{\bar{K}_j\hat{\bar{i}}(y_j)}\right|
               &\le& C_2,  \lb{4.17}\eea
where the sum in the left hand side of (\ref{4.16}) equals to
$\;N\sum_{\hat {\bar{i}}(y_j)>0}\frac{\hat{\bar{\chi}}(y_j)}{\hat{\bar{i}}(y_j)}\;$, the sum in the left hand side of
(\ref{4.17}) equals to $\;N\sum_{\hat{\bar{i}}(y_j)<0}\frac{\hat{\bar{\chi}}(y_j)}{\hat{\bar{i}}(y_j)}\;$ by (\ref{3.19}).}

In fact, we have the estimates
\bea  ^\#\{s\in\N_0 &|& 2C\le \bar{i}(y_j^{s\bar{K}_j+2m-1})+l\le 2N\}   \nn\\
&=& \;^\#\{s\in\N_0 \;| \; 2C\le \bar{i}(y_j^{s\bar{K}_j+2m-1})+l\le 2N,\;
             |\bar{i}(y_j^{s\bar{K}_j+2m-1})-(s\bar{K}_j+2m-1)\hat{\bar{i}}(y_j)|\le 2n\}  \nn\\
&\le& \;^\#\{s\in\N_0 \;| \;0< (s\bar{K}_j+2m-1)\hat{\bar{i}}(y_j)\le 2N-l+2n\}  \nn\\
&=& \;^\#\left\{s\in\N_0 \; \left |\;\frac{}{}\right.
     \;0\le s\le \frac{2N-l+2n-(2m-1)\hat{\bar{i}}(y_j)}{\bar{K}_j\hat{\bar{i}}(y_j)}\right\}  \nn\\
&\le& \; \frac{2N-l+2n}{\bar{K}_j\hat{\bar{i}}(y_j)}+1.  \nn\eea
On the other hand, we have
\bea
^\#\{s\in\N_0 &|& 2C\le \bar{i}(y_j^{s\bar{K}_j+2m-1})+l\le 2N\}  \nn\\
&=& \;^\#\{s\in\N_0 \;| \; 2C\le \bar{i}(y_j^{s\bar{K}_j+2m-1})+l\le 2N,\;
               |\bar{i}(y_j^{s\bar{K}_j+2m-1})-(s\bar{K}_j+2m-1)\hat{\bar{i}}(y_j)|\le 2n\}  \nn\\
&\ge& \;^\#\{s\in\N_0\;|\;\bar{i}(y_j^{s\bar{K}_j+2m-1})\le(s\bar{K}_j+2m-1)\hat{\bar{i}}(y_j)+2n\le 2N-l\} \nn\\
&\ge& \;^\#\left\{s\in\N_0 \; \left |\;\frac{}{}\right.
        \;0\le s\le \frac{2N-l-2n-(2m-1)\hat{\bar{i}}(y_j)}{\bar{K}_j\hat{\bar{i}}(y_j)}\right\}  \nn\\
&\ge& \;\frac{2N-l-2n}{\bar{K}_j\hat{\bar{i}}(y_j)}-2,  \nn\eea
where $m\le \bar{K}_j/2$ is used and we note that $\hat{\bar{i}}(y_j)> 0$ when
$2C\le \bar{i}(y_j^{s\bar{K}_j+2m-1})+l\le 2N$. Combining
these two estimates together with (\ref{4.14}), we obtain (\ref{4.16}). Similarly, we obtain (\ref{4.17}).

Note that all coefficients of $\bar{U}_a(t)$ in (\ref{4.8}) and (\ref{4.9}) are nonnegative. Hence, by Claim 1,
we have $\bar{u}_{h}\leq w_h\leq C_1$ for $h=2N$ or $-2C$ and $\bar{u}_{h}\leq w_{h+1}\leq C_1$ for $h=2C-1$ or $-2N-1$.
Now we choose $a$ to be sufficiently large, then we can choose $N$ to be sufficiently large.

Note that by Claims 1 and 2, the constants $C_1$ and $C_2$ are independent of $a$ and $K$. Hence
dividing both sides of (\ref{4.8}), (\ref{4.9}) by $N$ and letting $t=-1$, we obtain
\bea
\frac{\bar{M}_a^{2N}(2C;-1)-(N-C+1)}{N} &=& \frac{\bar{u}_{2N}+\bar{u}_{2C-1}}{N}, \nn\\
\frac{\bar{M}_a^{-2C}(-2N; -1)}{N} &=& \frac{\bar{u}_{-2N-1}+\bar{u}_{-2C}}{N}. \nn\eea
Let $N$ tend to infinity, then
\bea
&&\lim_{N\to\infty}\frac{1}{N}\bar{M}_a^{2N}(2C; -1) = 1, \nn\\
&&\lim_{N\to\infty}\frac{1}{N}\bar{M}_a^{-2C}(-2N; -1)= 0. \nn\eea
Hence (\ref{1.3}) and (\ref{1.4}) follow from (\ref{4.16}) and (\ref{4.17}).\hfill\hb

Let us also mention that if there is no solution with $\hat{\bar{i}}=0$, we do not need to cut our series at
$\pm 2C$; we can cut at $-2N$ and $2N$ only, thus obtaining
\be   \bar{M}(t)-\frac{1}{1-t^2}=(1+t)\bar{U}(t),   \lb{4.18}\ee
where $\bar{M}(t)$ denotes $\bar{M}_a(t)$ as $a$ tends to infinity, $\bar{U}(t)$ denotes $\bar{U}_a(t)$ as
$a$ tends to infinity.

\setcounter{equation}{0}
\section{A brief review on the mean index identities for closed characteristics on compact star-shaped
hypersurfaces in $\R^{2n}$}%{Section 5}

In this section, we briefly review the mean index identities for
closed characteristics on $\Sg\in\H_{st}(2n)$ developed in \cite{LLW1} which will be needed in Section 6. All
the details of proofs can be found in \cite{LLW1}.

We fix a $\Sg\in\H_{st}(2n)$ and suppose the condition (F) at the beginning of section 2 holds.
Note that by the condition (F), the number of critical orbits of $\Psi_{a, K}$
is finite. Hence as usual we can make the following definition.

{\bf Definition 5.1.} {\it Suppose $u$ is a nonzero critical point of $\Psi_{a, K}$, and $\Nn$ is an $S^1$-invariant open
neighborhood of $S^1\cdot u$ such that $\crit(\Psi_{a,K})\cap (\Lm_{a,K}(u)\cap \Nn) = S^1\cdot u$.
Then the $S^1$-critical
modules of $S^1\cdot u$ are defined by
\bea C_{S^1,\; q}(\Psi_{a, K}, \;S^1\cdot u)
=H_{q}((\Lambda_{a, K}(u)\cap\Nn)_{S^1},\; ((\Lambda_{a,K}(u)\setminus S^1\cdot u)\cap\Nn)_{S^1}). \nn\eea
Similarly, we define the $S^1$-critical modules $C_{S^1,\; q}(F_{a, K}, \;S^1\cdot x)$ of $S^1\cdot x$
for $F_{a, K}$.}

We fix $a$ and let $u_K\neq 0$ be a critical point of $\Psi_{a, K}$ with multiplicity $\mul(u_K)=m$,
that is, $u_K$ corresponds to a closed characteristic $(\sigma, z)\subset\Sigma$ with $(\sigma, z)$
being $m$-iteration of
some prime closed characteristic. Precisely, we have $u_K=-J\dot x+Kx$ with $x$
being a solution of (\ref{2.1}) and $x=\rho z(\frac{\sigma t}{T})$ with
$\frac{\vf_a^\prime(\rho)}{\rho}=\frac{\sigma}{aT}$.
Moreover, $(\sigma, z)$ is a closed characteristic on $\Sigma$ with minimal period $\frac{\sigma}{m}$.
By Lemma 2.10 of \cite{LLW1}, we construct a finite dimensional $S^1$-invariant subspace $G$ of
$L^{2}({\bf R}/{T {\bf Z}}; {\bf R}^{2n})$ and a functional $\psi_{a,K}$ on $G$.
For any $p\in\N$ satisfying $p\sigma<aT$, we choose $K$
such that $pK\notin \frac{2\pi}{T}\Z$, then the $p$th iteration $u_{pK}^p$ of $u_K$ is given by $-J\dot x^p+pKx^p$,
where $x^p$ is the unique solution of (\ref{2.1}) corresponding to $(p\sigma, z)$
and is a critical point of $F_{a, pK}$, that
is, $u_{pK}^p$ is the critical point of $\Psi_{a, pK}$ corresponding to $x^p$.
Denote by $g_{pK}^p$ the critical point of $\psi_{a,pK}$ corresponding to $u_{pK}^p$
and let $\wtd{\Lambda}_{a,K}(g_K)=\{g\in G \;|\; \psi_{a, K}(g)\le\psi_{a, K}(g_K)\}$.

Now we use the theory of Gromoll and Meyer, denote by $W(g_{pK}^p)$ the local characteristic manifold
of $g_{pK}^p$. Then we have

{\bf Proposition 5.2.}(cf. Proposition 4.2 of \cite{LLW1})
{\it For any $p\in\N$, we choose $K$ such that $pK\notin \frac{2\pi}{T}\Z$. Let $u_K\neq 0$
be a critical point of $\Psi_{a, K}$ with $\mul(u_K)=1$, $u_K=-J\dot x+Kx$ with $x$ being a critical point of
$F_{a, K}$. Then for all $q\in\Z$, we have
\bea
&& C_{S^1,\; q}(\Psi_{a,pK},\;S^1\cdot u_{pK}^p) \nn\\
&&\quad\cong \left(\frac{}{}H_{q-i_{pK}(u_{pK}^p)}(W(g_{pK}^p)\cap \wtd{\Lambda}_{a,pK}(g_{pK}^p),(W(g_{pK}^p)
                 \setminus\{g_{pK}^p\})\cap \wtd{\Lambda}_{a,pK}(g_{pK}^p))\right)^{\beta(x^p)\Z_p},
                  \qquad\quad \lb{5.1}\eea
where $\beta(x^p)=(-1)^{i_{pK}(u_{pK}^p)-i_K(u_K)}=(-1)^{i^v(x^p)-i^v(x)}$. Thus
\bea C_{S^1,\; q}(\Psi_{a,pK},\;S^1\cdot u_{pK}^p)=0 \quad
if q<i_{pK}(u_{pK}^p)~ or ~q>i_{pK}(u_{pK}^p)+\nu_{pK}(u_{pK}^p)-1.\lb{5.2}\eea
In particular, if $u_{pK}^p$ is
non-degenerate, i.e., $\nu_{pK}(u_{pK}^p)=1$, then}
\be C_{S^1,\; q}(\Psi_{a,pK},\;S^1\cdot u_{pK}^p)
    = \left\{\matrix{\Q, & {\rm if\;}q=i_{pK}(u_{pK}^p)\;{\rm and\;}\beta(x^p)=1,  \cr
                      0, & {\rm otherwise}. \cr}\right. \lb{5.3}\ee
We make the following definition:

{\bf Definition 5.3.} {\it For any $p\in\N$, we choose $K$ such that $pK\notin \frac{2\pi}{T}\Z$. Let $u_K\neq 0$ be
a critical point of $\Psi_{a,K}$ with $\mul(u_K)=1$, $u_K=-J\dot x+Kx$ with $x$ being a critical point of $F_{a, K}$.
Then for all $l\in\Z$, let
\bea
k_{l,\pm 1}(u_{pK}^p) &=& \dim\left(\frac{}{}H_{l}(W(g_{pK}^p)\cap
  \wtd{\Lambda}_{a,pK}(g_{pK}^p),(W(g_{pK}^p)\bs\{g_{pK}^p\})\cap \wtd{\Lambda}_{a,pK}(g_{pK}^p))\right)^{\pm\Z_p},
           \quad  \nn\\
k_l(u_{pK}^p) &=& \dim\left(\frac{}{}H_{l}(W(g_{pK}^p)\cap
  \wtd{\Lambda}_{a,pK}(g_{pK}^p),(W(g_{pK}^p)\bs\{g_{pK}^p\})\cap \wtd{\Lambda}_{a,pK}(g_{pK}^p))\right)^{\beta(x^p)\Z_p}.
           \qquad\quad   \nn\eea
Here $k_l(u_{pK}^p)$'s are called critical type numbers of $u_{pK}^p$. }

By Theorem 3.3 of \cite{LLW1}, we obtain that $k_l(u_{pK}^p)$ is independent of the choice of $K$ and
denote it by $k_l(x^p)$, here $k_l(x^p)$'s are called critical type numbers of $x^p$.

We have the following properties for critical type numbers:

{\bf Proposition 5.4.}(cf. Proposition 4.6 of \cite{LLW1})
{\it Let $x\neq 0$ be a critical point of $F_{a,K}$ with $\mul(x)=1$ corresponding to a
critical point $u_K$ of $\Psi_{a, K}$. Then there exists a minimal $K(x)\in \N$ such that
\bea
&& \nu^v(x^{p+K(x)})=\nu^v(x^p),\quad i^v(x^{p+K(x)})-i^v(x^p)\in 2\Z,  \qquad\forall p\in \N,  \lb{5.4}\\
&& k_l(x^{p+K(x)})=k_l(x^p), \qquad\forall p\in \N,\;l\in\Z. \lb{5.5}\eea
We call $K(x)$ the minimal period of critical modules of iterations of the functional $F_{a, K}$ at $x$. }

For every closed
characteristic $(\sigma,z)$ on $\Sigma$, let $aT>\sigma$ and choose $\vf_a$ as above.
Determine $\rho$ uniquely by $\frac{\vf_a'(\rho)}{\rho}=\frac{\sigma}{aT}$. Let $x=\rho z(\frac{\sigma t}{T})$.
Then we define the index $i(\sigma,z)$ and nullity $\nu(\sigma,z)$ of $(\sigma,z)$ by
$$ i(\sigma,z)=i^v(x), \qquad \nu(\sigma,z)=\nu^v(x). $$
Then the mean index of $(\sigma,z)$ is defined by
\be \hat i(\sigma,z) = \lim_{m\rightarrow\infty}\frac{i(m\sigma,z)}{m}.  \lb{5.6}\ee

Note that by Proposition 2.11 of \cite{LLW1}, the index and nullity are well defined and are independent of the
choice of $a$.

For a closed characteristic $(\sigma,z)$ on $\Sigma$, we simply denote by $z^m\equiv(m\sigma,z)$
the m-th iteration of $z$ for $m\in\N$.
By Proposition 3.2 of \cite{LLW1}, we can define the critical type numbers $k_l(z^m)$ of $z^m$ to be $k_l(x^m)$,
where $x^m$ is the critical point of $F_{a, K}$ corresponding to $z^m$. We also define $K(z)=K(x)$.

{\bf Lemma 5.5.} {\it For a symmetric closed characteristic $(\sigma,z)$,
we have $\hat{\bar{i}}(z)=\frac{\hat{i}(z)}{2}$, where $\hat{\bar{i}}(z)$
is defined in Definition 3.6.}

{\bf Proof.} Denote by $\gamma\equiv
\gamma_{z}$ the associated symplectic paths of $(\sigma,z)$, let $\psi=\gamma|_{[0, \sigma/2]}$, then we have
$\gamma|_{[0, \sigma]}=\psi^2$. Now using Theorem 9.2.1 of \cite{Lon3} and Definition 3.6,
we obtain \bea   \hat{\bar{i}}(z)
&=& \lim_{m\to\infty}\frac{\bar{i}(z^{2m-1})}{2m-1}    \nn\\
&=& \lim_{m\to\infty}\frac{i_{-1}(z^{2m-1})}{2m-1}    \nn\\
&=& \lim_{m\to\infty}\frac{i_1(\psi^{4m-2})-i_1(\psi^{2m-1})}{2m-1}   \nn\\
&=& \hat{i}(\psi)\;=\;\frac{\hat{i}(z)}{2},\nn\eea
our lemma follows. \hfill\hb

{\bf Remark 5.6.}(cf. Remark 4.10 of \cite{LLW1}) {\it
Note that $k_l(z^m)=0$ for $l\notin [0, \nu(z^m)-1]$ and it can take only values $0$
or $1$ when $l=0$ or $l=\nu(z^m)-1$. Moreover, the following facts are useful:

(i) $k_0(z^m)=1$ implies $k_l(z^m)=0$ for $1\le l\le \nu(z^m)-1$.

(ii) $k_{\nu(z^m)-1}(z^m)=1$ implies $k_l(z^m)=0$ for $0\le l\le \nu(z^m)-2$.

(iii) $k_l(z^m)\ge 1$ for some $1\le l\le \nu(z^m)-2$ implies $k_0(z^m)=k_{\nu(z^m)-1}(z^m)=0$.

(iv) In particular, only one of the $k_l(z^m)$s for $0\le l\le \nu(z^m)-1$ can be non-zero when $\nu(z^m)\le 3$.}

For a closed characteristic $(\sigma,z)$ on $\Sigma$, the average Euler characteristic $\hat\chi(z)$ of $z$ is
defined by\bea \hat\chi(z)=\frac{1}{K(z)}\sum_{1\le m\le K(z)\atop 0\le l\le 2n-2}
(-1)^{i(z^{m})+l}k_l(z^{m}).  \lb{5.7}\eea
$\hat\chi(z)$ is a rational number. In particular, if all $z^m$s are
non-degenerate, then by Proposition 5.4 we have
\be \hat\chi(z)
    = \left\{\matrix{(-1)^{i(z)}, & {\rm if\;\;} i(z^2)-i(z)\in 2\Z,  \cr
           \frac{(-1)^{i(z)}}{2}, & {\rm otherwise}. \cr}\right.  \lb{5.8}\ee
We have the following mean index identities for closed characteristics.

{\bf Theorem 5.7.} {\it Suppose that $\Sg\in \H_{st}(2n)$ satisfies $\,^{\#}\T(\Sg)<+\infty$. Denote all
the geometrically distinct prime closed characteristics by $\{(\sigma_j,\; z_j)\}_{1\le j\le k^\prime}$. Then the
following identities hold
\bea
\sum_{1\le j\le k^\prime\atop \hat{i}(z_j)>0}\frac{\hat{\chi}(z_j)}{\hat{i}(z_j)} &=& \frac{1}{2}, \lb{5.9} \\
\sum_{1\le j\le k^\prime\atop \hat{i}(z_j)<0}\frac{\hat{\chi}(z_j)}{\hat{i}(z_j)} &=& 0. \lb{5.10}\eea}

Let $F_{a, K}$ be a functional defined by (\ref{2.4}) for some $a, K\in\R$ sufficiently large and let $\epsilon>0$ be
small enough such that $[-\epsilon, 0)$ contains no critical values of $F_{a, K}$. For $b$ large enough,
The normalized Morse series of $F_{a, K}$ in $ X^{-\ep}\setminus X^{-b}$
is defined, as usual, by
\be  M_a(t)=\sum_{q\ge 0,\;1\le j\le p} \dim C_{S^1,\;q}(F_{a, K}, \;S^1\cdot v_j)t^{q-d(K)},  \lb{5.11}\ee
where we denote by $\{S^1\cdot v_1, \ldots, S^1\cdot v_p\}$ the critical orbits of $F_{a, K}$ with critical
values less than $-\epsilon$. The Poincar\'e series of
$H_{S^1, *}( X, X^{-\ep})$ is $t^{d(K)}Q_a(t)$, according to Theorem 5.1 of \cite{LLW1}, if we set
$Q_a(t)=\sum_{k\in \Z}{q_kt^k}$, then
\be   q_k=0 \qquad\qquad \forall\;k\in \mathring {I},  \lb{5.12}\ee
where $I$ is an interval of $\Z$ such that $I \cap [i(\sigma, z), i(\sigma, z)+\nu(\sigma, z)-1]=\emptyset$ for all
closed characteristics $(\sigma, z)$ on $\Sigma$ with $\sigma\ge aT$. Then by Section 6 of \cite{LLW1}, we have
\be  M_a(t)-\frac{1}{1-t^2}+Q_a(t) = (1+t)U_a(t),   \lb{5.13}\ee
where $U_a(t)=\sum_{i\in \Z}{u_it^i}$ is a Laurent series with nonnegative coefficients.
If there is no closed characteristic with $\hat{i}=0$, then \be   M(t)-\frac{1}{1-t^2}=(1+t)U(t),   \lb{5.14}\ee
where $M(t)=\sum_{i\in\Z}{m_it^i}$ denotes the limit of $M_a(t)$ as $a$ tends to infinity,
$U(t)=\sum_{i\in\Z}{u_it^i}$ denotes the limit of $U_a(t)$ as $a$ tends to infinity and possesses
only non-negative coefficients. Specially for our later applications in Section 6, suppose that there exists an
integer $p<0$ such that the coefficients of $M(t)$ satisfy $m_p>0$ and $m_q=0$ for all integers $q<p$. Then
(\ref{5.14}) implies
\be   m_{p+1} \ge m_p.  \lb{5.15}\ee

\setcounter{equation}{0}%\setcounter{figure}{0}
\section{Proof of Theorem 1.4}%{Section 6}

In this section, we prove Theorem 1.4 by using the mean index identities in Theorem 1.1 and Theorem 5.7, Morse
inequality and the index iteration theory developed by Y. Long and his coworkers.

The following theorem relates the Morse index defined in Section 5 to the Maslov-type index.

{\bf Theorem 6.1.} (cf. Theorem 2.1 of \cite{HuL1}) {\it Suppose $\Sg\in \H_{st}(2n)$ and
$(\tau,y)\in \T(\Sigma)$. Then we have
\be i(y^m)\equiv i(m\tau,y)=i(y, m)-n,\quad \nu(y^m)\equiv\nu(m\tau, y)=\nu(y, m),
       \qquad \forall m\in\N, \lb{6.1}\ee
where $i(y, m)$ and $\nu(y, m)$ are the Maslov-type index and nullity of $(m\tau,y)$ (cf. Section 5.4 of \cite{Lon3}).
In particular, we have $\hat{i}(\tau,y)=\hat{i}(y,1)$, where $\hat{i}(\tau ,y)$ is given in Section 5, $\hat{i}(y,1)$
is the mean Maslov-type index (cf. Definition 8.1 of \cite{Lon3}). Hence we denote it simply by $\hat{i}(y)$.}

In the following, we fix $n=2$. Before we give the proof of Theorem 1.4, we need a proposition:

{\bf Proposition 6.2.} {\it Let $\Sg\in \H_{st}(4)$ satisfy $\,^{\#}\T(\Sg)=2$. Denote
the two geometrically distinct prime closed characteristics by $\{(\tau_j,\; y_j)\}_{1\le j\le 2}$. If
$i(y_j)\geq 0, j=1, 2$, then both of the closed characteristics are elliptic.}

{\bf Proof.} Denote by $\gamma_j\equiv \gamma_{y_j}$ the associated symplectic paths of
$(\tau_j,\; y_j)$ for $1\le j\le 2$. Then by Lemma 3.3 of \cite{HuL1} (cf. also Lemma 15.2.4 of
\cite{Lon3}), there exist $P_j\in \Sp(4)$ and $M_j\in \Sp(2)$ such that
\be \gamma_j(\tau_j)=P_j^{-1}(N_1(1,\,1)\diamond M_j)P_j,~j=1, 2, \lb{6.2}\ee
where $N_1(\lambda,\,b)=\left(
                    \begin{array}{cc}
                      \lambda & b\\
                      0 & \lambda \\
                    \end{array}\right)$ for $\lambda, b\in\R$.
By assumption that $i(y_j)\geq 0$ and Theorem 6.1, we have $i(y_j, 1)\geq 2$, together with (\ref{6.2}),
by Corollary 8.3.2 of \cite{Lon3}, it gives $\hat{i}(y_j)>2$. Then (\ref{5.14})
holds. Thus for every $k\in \N$, when $a$ is large enough, there exist some $1\leq j \leq 2$ and $m\in\N$ such that
\be  C_{S^1,\;d(K)+2(k-1)}(F_{a,K},\;S^1\cdot {x}_j^m)\neq 0, \lb{6.3}\ee
where $x_j$ is the critical point of $F_{a,K}$ corresponding to $y_j$.
Using the common index jump theorem (Theorems 4.3 and 4.4 of \cite{LoZ1}, Theorems 11.2.1
and 11.2.2 of \cite{Lon3}), we obtain infinitely many $(N, m_1, m_2)\in\N^3$ such that
\bea i(y_j, 2m_j+1) &=& 2N+i(y_j, 1),\lb{6.4}\\
i(y_j, 2m_j-1)+\nu(y_j, 2m_j-1) &=& 2N-(i(y_j, 1)+2S^+_{\gamma_j(\tau_j)}(1)-\nu(y_j, 1)).\lb{6.5}\eea
Since $i(y_j, 1)\geq 2$, then (22) on Page 340 of \cite{Lon3} also holds for $n=2$ and (\ref{6.4}),
(\ref{6.5}) yield
\bea i(y_j, 2m_j+1) &\geq & 2N+2,\lb{6.6}\\
i(y_j, 2m_j-1)+\nu(y_j, 2m_j-1)-1 &\leq & 2N-3.\lb{6.7}\eea
Combining $i(y_j, 1)\geq 2$ with Theorem 10.2.4 of \cite{Lon3}, we have
\be  i(y_j, m)<i(y_j, m+1), \qquad \forall\;m\in\N,\; j=1,2. \lb{6.8}\ee
Thus by (\ref{6.3}), (\ref{6.6})-(\ref{6.8}) and Theorem 6.1, noticing Remark 5.6 (iv), for $s=1, 2$
we have $C_{S^1,\;2N-2s}(F_{a,K},\;S^1\cdot {x}_{j_s}^{2m_{j_s}})\neq 0$, where
$\{j_s\mid s=1, 2\}=\{1, 2\}$. Then by the proof of Theorem 15.5.2 of \cite{Lon3},
we obtain that the two closed characteristics are elliptic.\hfill\hb

{\bf Remark 6.3.} {\it If $\Sg\in \H_{con}(4)$ satisfy $\,^{\#}\T(\Sg)=2$, then $i(y_j)\geq 0, j=1, 2$ holds.
Thus both of the closed characteristics are elliptic by the above proposition. But for star-shaped case,
such a good lower bound for initial indices may not hold generally. Thus in the following proof of Theorem 1.4,
we continue in four cases and try to get such a lower bound for all cases.}

{\bf Proof of Theorem 1.4.} Let $\Sg \in \mathcal {S}\H_{st}(4)$ and denote by $\{(\tau_1,\; y_1), \;(\tau_2,\; y_2)\}$
the two geometrically distinct prime closed characteristics on $\Sg $, and by $\ga_j\equiv \ga_{y_j}$ the associated
symplectic paths of $(\tau_j,\; y_j)$ for $1\le j\le 2$. Then as in the proof of Proposition 6.2, there exist
$P_j\in \Sp(4)$ and $M_j\in \Sp(2)$ such that
\be   \ga_j(\tau_j)=P_j^{-1}(N_1(1,\,1)\diamond M_j)P_j, \quad {\rm for}\;\;j=1, 2.  \lb{6.9}\ee
Since $\Sg \in \mathcal{S}\H_{st}(4)$, by Theorem 1 of \cite{Gir1} and Lemma 4.2 of \cite{LLZ1}, both $(\tau_1,\; y_1)$
and $(\tau_2,\; y_2)$ must be symmetric. Let $\psi_j=\gamma_j|_{[0, \tau_j/2]}$, then we have
\be   \gamma_j|_{[0, \tau_j]}=\psi_j^2 \quad {\rm and}\quad \gamma_j(\tau_j)=\gamma_j(\tau_j/2)^2. \lb{6.10}\ee

Note that by Section 9 of \cite{Vit2}, we know that there exists at least one non-hyperbolic closed characteristic
on $\Sg$ and it is certainly elliptic when $n=2$. In the following, we prove Theorem 1.4 by contradiction. Without
loss of generality, we suppose that $(\tau_1, y_1)$ is elliptic and $(\tau_2, y_2)$ is hyperbolic.

For the properties of these two orbits, we have

{\bf Claim 1.} {\it The closed characteristics $(\tau_1, y_1)$ and $(\tau_2, y_2)$ satisfy

(i) $i(y_2^m) = m(i(y_2)+3) - 3$ and $\nu(y_2^m)=1$ for all $m\in\N$, and thus $\hat{i}(y_2) = i(y_2) + 3$.

(ii) $\hat{i}(y_2) >0$;

(iii) $i(\psi_2^m) = m i(\psi_2) - \frac{1+(-1)^m}{2}$ for all $m\in\N$;

(iv) $i(y_2) = i(y_2,1) - 2 = i(\psi_2^2) - 2 = 2 i(\psi_2) - 3 \;\in\; 2\N-3$;

(v) $\hat{\chi}(y_2)=(-1)^{i(y_2)}=-1$ and
$$ \frac{\hat{\chi}(y_2)}{\hat{i}(y_2)}=\frac{-1}{\hat{i}(y_2)}< 0.  $$

(vi) If $i(y_2)=-1$, then $\{y_2^m\,|\,m\in\N\}$ contributes to every Morse type number $m_{q}$ in (\ref{5.14})
precisely a $1$ for each odd $q\in\Z$ with $q\ge i(y_2)$ and nothing for all the other $q\in\Z$.

(vii) $0< \hat{i}(y_1) \in \Q$.

(viii) $\hat{\chi}(y_1)>0$ and
$$  \frac{\hat{\chi}(y_1)}{\hat{i}(y_1)} = \frac{1}{2} - \frac{\hat{\chi}(y_2)}{\hat{i}(y_2)} > \frac{1}{2}.  $$}

In fact, by Theorem 8.3.1 of \cite{Lon3} (Theorem 1.3 of \cite{Lon2}), we have $i(y_2,m) = m(i(y_2, 1)+1)-1$
for all $m\in\N$. Together with Theorem 6.1, we obtain (i).

We claim $\hat{i}(y_2)\neq 0$. In fact, because $y_2$ is hyperbolic, $y_2^m$ is non-degenerate for every
$m\ge 1$. Thus if $\hat{i}(y_2)=0$, we then have $i(y_2^m)=i(y_2,m)-2=-3$ for all $m\ge 1$. Then the Morse
type number satisfies $M_{-3}=+\infty$. But then $\hat{i}(y_1)$ must be positive by Theorem 5.7, and
contributions of $\{y_1^m\}$ to every Morse type number thus must be finite. Then the Morse inequality yields a
contradiction and proves the claim (cf. the proof below (9.3) of \cite{Vit2} for details).

If $\hat{i}(y_2)<0$, by (\ref{5.10}) we obtain
\be    \frac{\hat{\chi}(y_2)}{\hat{i}(y_2)} = 0.   \lb{6.11}\ee
But because $(\tau_2, y_2)$ is hyperbolic, by (\ref{5.8}) we have $\hat{\chi}(y_2)\not= 0$, which
contradicts to (\ref{6.11}) and proves (ii).

Because $(\tau_2, y_2)$ is hyperbolic, the $2\times 2$ matrix $M_2$ is hyperbolic. Thus by the proof of
Proposition 2.12, we have
\be   \ga_2(\tau_2/2) = N_1(-1,-1)\dm C,  \lb{6.12}\ee
in an appropriate coordinates and the fact $C^2=M_2$ implies $\sg(C)\cap \U = \emptyset$, where $C$ is
a symplectic matrix and $\sg(C)$ denotes the spectrum of the matrix $C$, ${\bf U}$ is the unit circle
in the complex plane. Then by Theorem 8.3.1 of \cite{Lon3}, we obtain (iii) and then $i(y_2)$ is odd. Then
together with (i) and (ii), it yields $i(y_2)\ge -1$, i.e., (iv) holds.

Since $i(y_2^2)-i(y_2)=i(y_2)+3\in 2\Z$ holds and $y_2$ is hyperbolic by the above (i) and (iv), by
(\ref{5.8}) and the above (iv) we obtain the first equality in (v). Together with above (ii), it yields
the second equality and estimate in (v).

Note that if $i(y_2)=-1$, by (i), the set $\Th(y_2)=\{i(y_2^m)\,|\,m\in\N\}$ consists of every odd integer not
less than $i(y_2)$ precisely once. Note also that $y_2^m$ is non-degenerate for every $m\in\N$ by (i). Then (vi)
follows from (\ref{5.3}) of Proposition 5.2.

If $(\tau_1, y_1)$ and its iterates are all non-degenerate, since $(\tau_1, y_1)$ is elliptic, then
$\hat{i}(y_1)$ must be irrational by Corollary 8.3.2 of \cite{Lon3} and then so is
$\frac{\hat{\chi}(y_1)}{\hat{i}(y_1)}$, because $\hat{\chi}(y_1)$ is rational. Then by (\ref{5.9})
of Theorem 5.7, the other closed characteristic $(\tau_2, y_2)$ must possess an irrational mean index
$\hat{i}(y_2)$, which contradicts to the second identity in (i), and thus $\hat{i}(y_1)$ must be
rational. On the other hand, because of the above (ii) and (v) as well as (\ref{5.9}) of Theorem 5.7,
we obtain $\hat{i}(y_1)>0$ and proves (vii).

Now by the above (ii), (v), (vii), and (\ref{5.9}) of Theorem 5.7 we obtain (viii).

The proof of Claim 1 is complete.

By (vii) of Claim 1, we only need to consider the following four cases according to the classification
of basic norm forms of $\ga_1(\tau_1)$. In the following we use the notations from Definition 1.8.5
and Theorem 1.8.10 of \cite{Lon3}, and specially we let $R(\th) =  \left(\begin{array}{cc}
                                                             \cos{\th} & -\sin{\th} \\
                                                             \sin{\th} &  \cos{\th} \\
                                                                \end{array}\right)$ with $\th\in\R$,
and use $M\dm N$ to denote the symplectic direct sum of two symplectic matrices $M$ and $N$ as
in pages 16-17 of \cite{Lon3}.

\vskip 2 mm

{\bf Case 1.} {\it $\ga_1(\tau_1)$ can be connected to $N_1(1,1)\dm N_1(-1,b)$ within
$\Om^0(\ga_1(\tau_1))$ with $b=0$ or $\pm 1$.}

In this case, by Theorems 8.1.4 and 8.1.5 of \cite{Lon3}, and Theorem 6.1, we have
\be  i(y_1, 1) \quad {\rm and}\quad i(y_1) \quad {\rm are\;even}.  \lb{6.13}\ee
By Theorem 1.3 of \cite{Lon2}, we have
\bea
i(y_1,m) &=& m(i(y_1, 1)+1)-1, \quad {\rm for}\;\; b=1;\nn\\
i(y_1, m) &=& m(i(y_1,1)+1)-1-\frac{1+(-1)^m}{2}, \quad {\rm for}\;\;b=0, -1. \nn\eea
By Theorem 6.1, we obtain
\bea
i(y_1^m) &=& m(i(y_1)+3)-3, \quad {\rm for}\;\; b=1;   \lb{6.14}\\
i(y_1^m) &=& m(i(y_1)+3)-3-\frac{1+(-1)^m}{2}, \quad {\rm for}\;\;b=0, -1. \lb{6.15}\eea
Then in both cases we obtain
\be \hat{i}(y_1) = i(y_1)+3.  \lb{6.16}\ee

Note that by Proposition 5.4 and the form of $\ga_1(\tau_1)$, we have $K(y_1)=2$. Thus by (viii) of
Claim 1 and (\ref{5.7}), we obtain
\be   0 < \hat{\chi}(y_1)=\frac{1+(-1)^{i(y_1^2)}(k_0(y_1^2)-k_1(y_1^2)+k_2(y_1^2))}{2}.   \lb{6.17}\ee
Because at most one of $k_l(y_1^2)'s$ for $0\leq l\leq 2$ can be non-zero by Remark 5.6 (iv), we obtain
\be   (-1)^{i(y_1^2)+l}k_l(y_1^2)\ge 0, \qquad {\rm for}\;\;l=0, 1, 2.   \lb{6.18}\ee

when $m$ is odd, we have $\nu(x_1^m)=1$ by the assumption on $\ga_1(\tau_1)$. In this case, because
$i(y_1)$ is even by (\ref{6.13}), we have $i^v(x_1^m) = i(y_1^m) = m(i(y_1)+3)-3$ is even, and then
$$  \bb(x^m) = (-1)^{i^v(x_1^m)-i^v(x_1)} = 1, $$
where we denote by $x_j$ the critical point of $F_{a,K}$ corresponding to $y_j$ for $j=1$ and $2$.
Thus by (\ref{5.3}) of Proposition 5.2 for every odd $m\in\N$, we obtain
\bea
C_{S^1,\;d(K)+k}(F_{a,K},\;S^1\cdot {x}_1^{m})=\Q, && {\rm if}\;\;k=i(y_1^m),   \lb{6.19}\\
C_{S^1,\;d(K)+k}(F_{a,K},\;S^1\cdot {x}_1^{m})=0, && {\rm if}\;k\not=i(y_1^m),   \lb{6.20}\eea
where (\ref{6.20}) holds specially when $k\in 2\Z-1$.

When $m$ is even, we continue in two cases (A) for $b=1$ with (\ref{6.14}) and (B) for
$b=0, -1$ with (\ref{6.15}).

(A) {\it $m$ is even, $b=0$ or $-1$, and (\ref{6.15}) holds.}

In this case, $i(y_1^2)$ is even by (\ref{6.15}). Therefore by (\ref{6.17}) we obtain
\be  k_1(y^2) = 0, \quad \hat{\chi}(y_1)=\frac{1+(k_0(y_1^2)+k_2(y_1^2))}{2}>0.  \lb{6.21}\ee
Because $K(y_1)=2$, we then obtain
\be C_{S^1,\;d(K)+2k-1}(F_{a,K},\;S^1\cdot {x}_1^{m})=0, \qquad\forall k\in\Z, m\in 2\N.   \lb{6.22}\ee

Therefore when $b=0, -1$ from (\ref{6.19}), (\ref{6.20}) and (\ref{6.22}) we obtain
\be C_{S^1,\;d(K)+2k-1}(F_{a,K},\;S^1\cdot {x}_1^{m}) = 0, \qquad \forall k\in\Z, m\in \N. \lb{6.23}\ee

(B) {\it $m$ is even, $b=1$, and (\ref{6.14}) holds.}

In this case, $i(y_1^2)$ is odd by (\ref{6.14}). Therefore by (\ref{6.17}) we obtain
\be  k_0(y^2) = k_2(y^2) = 0, \quad 0<\hat{\chi}(y_1)=\frac{1+k_1(y_1^2)}{2}.  \lb{6.24}\ee
Because $K(y_1)=2$, we then obtain
\be C_{S^1,\;d(K)+2k-1}(F_{a,K},\;S^1\cdot {x}_1^{m})=0, \qquad \forall k\in\Z, m\in 2\N.   \lb{6.25}\ee
By (vii) and (viii) of Claim 1, we have
$$  \frac{1+k_1(y_1^2)}{2(i(y_1)+3)}
   = \frac{1}{2} + \frac{|\hat{\chi}(y_2)|}{i(y_2)+3} > \frac{1}{2}. $$
Specially this implies $k_1(y_1^2)>0$, and then when $m$ is even, we obtain
\be C_{S^1,\;d(K)+k}(F_{a,K},\;S^1\cdot {x}_1^{m})=\Q, \quad {\rm if\;and\;only\;if}\quad k=i(y_1^m)+1. \lb{6.26}\ee

Therefore when $b=1$, from (\ref{6.19}), (\ref{6.20}), (\ref{6.25}) and (\ref{6.26}), we obtain
\be C_{S^1,\;d(K)+2k-1}(F_{a,K},\;S^1\cdot {x}_1^{m}) = 0, \qquad \forall k\in\Z, m\in \N. \lb{6.27}\ee

Specially from (\ref{6.23}) and (\ref{6.27}), for any case we have
\be C_{S^1,\;d(K)+2k-1}(F_{a,K},\;S^1\cdot {x}_1^{m}) = 0, \qquad \forall k\in\Z, m\in \N. \lb{6.28}\ee

By (iv) of Claim 1, we have $i(y_2)\ge -1$ and it is odd. Then we continue our proof in two subcases:

{\bf Subcase 1.1.} $i(y_2)\ge 1$.

In this case, we have $i(y_2^m)\ge 0$ by (i) of Claim 1. Thus by (\ref{5.2}) of Proposition
5.2, we have $C_{S^1,\;d(K)-1}(F_{a,K},\;S^1\cdot {x}_2^{m})=0$ for all $m\in\N$. Combining it
with (\ref{6.28}), we get
\be   m_{-1}=0.   \lb{6.29}\ee
Here and below in this Section $m_i$ is the coefficient of $t^i$ of $M(t)=\sum_{i\in \Z}{m_it^i}$ in
(\ref{5.14}).

By (\ref{6.13}), (\ref{6.16}) and (vii) of Claim 1, we have $i(y_1)=-2$ or $i(y_1)\geq 0$.

If $i(y_1)=-2$, by Proposition 5.2, we have $C_{S^1,\;d(K)-2}(F_{a,K},\;S^1\cdot {x}_1)=\Q$.
Thus we have
\be    m_{-2}\ge 1.   \lb{6.30}\ee
Since $i(y_2)\ge 1$ and $i(y_1)=-2$, we get $C_{S^1,\;d(K)-q}(F_{a,K},\;S^1\cdot {x}_j^{m})=0$
for all $m\in\N$ and $q\geq 3$ with $j=1, 2$. Thus we obtain
\be    m_{-q}=0, \qquad \forall q\ge 3.   \lb{6.31}\ee
Then we have $m_{-1}\ge m_{-2}\ge 1$ by (\ref{6.30}), (\ref{6.31}) and (\ref{5.15}), which
contradicts to (\ref{6.29}).

If $i(y_1)\ge 0$. Noticing that $i(y_2)\ge 1$, from Proposition 6.2 we know that the two closed
characteristics are elliptic which contradicts to our assumption too.

{\bf Subcase 1.2.} $i(y_2)=-1$.

In this case, $\hat{i}(y_2)=2$ by (i) of Claim 1, and $i(\psi_2)=1$ by (iv) of Claim 1. Using
the Bott-type formulae (Theorem 9.2.1 of \cite{Lon3}) and (iii) of Claim 1, we have
\be  i_{-1}(\psi_2^m) = i(\psi_2^{2m})-i(\psi_2^m), \qquad \forall\;m\in\N.  \lb{6.32}\ee
Specially we obtain
$$   i_{-1}(\psi_2) = i(\psi_2^2)-i(\psi_2) = i(y_2) + 2 - i(\psi_2) = -1+2-1 = 0.  $$
From Definition 3.6, we have $\bar{i}(y_2^m)=i_{-1}(\psi_2^m)$ and then from (\ref{6.32}) we obtain
\be   \bar{i}(y_2)=i_{-1}(\psi_2)=0.   \lb{6.33}\ee
By Lemma 5.5, we have
\bea   \hat{\bar{i}}(y_2)
&=&\frac{\hat{i}(y_2)}{2} \;=\; 1.   \lb{6.34}\\
\hat{\bar{i}}(y_1)&=&\frac{\hat{i}(y_1)}{2}>0. \lb{6.35}\eea

Then we can apply Theorem 1.1 to get
\be  \frac{\hat{\bar{\chi}}(y_1)}{\hat{\bar{i}}(y_1)}
     + \frac{\hat{\bar{\chi}}(y_2)}{\hat{\bar{i}}(y_2)} = 1.  \lb{6.36}\ee
By (\ref{3.20}) and (\ref{6.33}), we have
\be  \hat{\bar{\chi}}(y_2)=(-1)^{\bar{i}(y_2)}=1.  \lb{6.37}\ee
It then follows from (\ref{6.34})-(\ref{6.37}) that
\be   \hat{\bar{\chi}}(y_1)=0.  \lb{6.38}\ee
Because in this case $\bar{\nu}(y_1^m)=1$ holds for all $m\in 2\N-1$, thus we can use
(\ref{3.20}) again, and obtain $\hat{\bar{\chi}}(y_1)=(-1)^{\bar{i}(y_1)}\neq 0$,
which contradicts to (\ref{6.38}). This completes the proof of Case 1.

\vskip 2 mm

{\bf Case 2.} $\ga_1(\tau_1)$ can be connected to $N_1(1,1)\dm R(\th)$ within $\Om^0(\ga_1(\tau_1))$
with some $\th\in (0,\pi)\cup (\pi,2\pi)$ and $\th/\pi\in \Q$.

In this case, we have always $K(y_1)\ge 3$ by the definition of $\th$. By Theorems 8.1.4 and 8.1.7
of \cite{Lon3} and Theorem 6.1 we obtain
\be  i(y_1,1)\quad {\rm and}\quad i(y_1) \quad {\rm are\;even}.  \lb{6.39}\ee
By Theorem 1.3 of \cite{Lon2} (i.e., Theorem 8.3.1 of \cite{Lon3}), we have
$$   i(y_1,m) = mi(y_1,1) + 2E(\frac{m\th}{2\pi}) - 2.  $$
By Theorem 6.1, we obtain
\be  i(y_1^m) = m(i(y_1)+2) + 2E(\frac{m\th}{2\pi}) - 4.   \lb{6.40}\ee
Then
\be   \hat{i}(y_1) = i(y_1) + 2 + \frac{\th}{\pi}.   \lb{6.41}\ee
In this case, by (ii) and (v)-(viii) of Claim 1 we obtain
\be  \frac{\hat{\chi}(y_1)}{\hat{i}(y_1)} = \frac{1}{2} + \frac{1}{\hat{i}(y_2)} > \frac{1}{2}.
             \lb{6.42}\ee
On the other hand, by (\ref{5.7}) we have
\be  \hat{\chi}(y_1)
\;=\; \frac{K(y_1)-1+k_0(y_1^{K(y_1)})-k_1(y_1^{K(y_1)})+k_2(y_1^{K(y_1)})}{K(y_1)} \;\le\; 1.
             \lb{6.43}\ee
Thus we obtain
\be   0<\hat{i}(y_1)<2.  \lb{6.44}\ee
Together with (\ref{6.39}) and (\ref{6.41}), it yields
\be  i(y_1) = -2.   \lb{6.45}\ee
Then by (\ref{6.40}) we obtain that
\bea
&& i(y_1^m) = 2E(\frac{m\th}{2\pi}) - 4 \ge -2, \quad \forall\;m\in\N, \lb{6.46}\\
&& \nu(y_1^m) = 1, \quad {\rm if}\;m\not= 0\;\;{\rm mod}\;\;K(y_1),  \lb{6.47}\\
&& \nu(y_1^{K(y_1)}) = 3.  \lb{6.48}\eea
Specially we obtain $\frac{K(y_1)\th}{2\pi}\in\N$ and
\be  \hat{i}(y_1)=\frac{\theta}{\pi}\ge \frac{2}{K(y_1)}. \lb{6.49}\ee

Because $\nu(y_1)=1$, by Proposition 5.2 we have $C_{S^1,\;d(K)-2}(F_{a,K},\;S^1\cdot {x}_1)=\Q$.
This proves
\be   m_{-2}\ge 1.  \lb{6.50}\ee
Note that by (i) and (ii) of Claim 1, we have $i(y_2)\ge -1$, and then $i(y_2^m)\ge -1$ for all
$m\in\N$. Thus we have
$$  C_{S^1,\;d(K)-q}(F_{a,K},\;S^1\cdot {x}_j^{m})=0, \qquad \forall\;m\in\N, \;q\ge 3,\; j=1,2.  $$
Thus we have
\be  m_{-q}=0, \qquad \forall\;q\ge 3.  \lb{6.51}\ee
By (\ref{5.15}), we get
\be   m_{-1}\ge 1.   \lb{6.52}\ee

In order to make (\ref{6.52}) hold, there are two possibilities. The first is that $y_1$
contributes a positive integer to $m_{-1}$, by (\ref{6.46}) and Proposition 5.2 which needs
\be \left\{\begin{array}{cc}
&i(y_1^m)=-2, \qquad \forall\;1\le m\le K(y_1), \\
&k_1(y_1^{K(y_1)}) > 0 \quad {\rm and}\quad k_0(y_1^{K(y_1)})=k_2(y_1^{K(y_1)})=0.  \\
                  \end{array}\right.    \lb{6.53}\ee
The second is that $y_2$ contributes a $1$ to $m_{-1}$, which requires $i(y_2)= -1$ by (iv)
of Claim 1.

In the second possibility of $i(y_2)= -1$, by (i) of Claim 1 we have $\hat{i}(y_2)=2$. Then
by (\ref{6.42}) and (\ref{6.43}), we obtain $\hat{i}(y_1)=\hat{\chi}(y_1)\le 1$. Because
$\th\neq\pi$, together with (\ref{6.41}) and (\ref{6.45}), it implies
$\frac{\th}{\pi} = \hat{i}(y_1) < 1$, which implies $i(y_1^2)=-2$ by (\ref{6.46}).
By (\ref{6.47}) and Proposition 5.2, we then obtain $C_{S^1,\;d(K)-2}(F_{a,K},\;S^1\cdot {x}_1^m)=\Q$
for $m=1$ and $2$, which implies $m_{-2}\ge 2$, and then $m_{-1}\ge 2$ by (\ref{5.15}) and
(\ref{6.51}). Now by the fact $i(y_2)= -1$ and (vi) of Claim 1, $\{y_1^m\}$ needs
to contribute to $m_{-1}$ too, and then we must have
$C_{S^1,\;d(K)-1}(F_{a,K},\;S^1\cdot {x}_1^m)\neq 0$ for some $m\in\N$. This implies that
(\ref{6.53}) holds always.

Then (\ref{6.53}) implies $K(y_1)\th = 2\pi$, and by (\ref{6.46}) we obtain
\be  \hat{i}(y_1) = \frac{\th}{\pi} = \frac{2}{K(y_1)}.  \lb{6.54}\ee
Therefore by (\ref{6.47}), (\ref{6.53}), (vi) of Claim 1, Proposition 5.2, and (\ref{5.15}),
noticing that $y_2^m$ contribute at most a 1 to $m_{-1}$, we have
$$   k_1(y_1^{K(y_1)})+1 \ge m_{-1} \ge m_{-2} = K(y_1)-1. $$
Therefore by (\ref{6.43}) we obtain
$$  \hat{\chi}(y_1) \le \frac{1}{K(y_1)}. $$
Together with (\ref{6.49}) we then obtain
$$  \frac{\hat{\chi}(y_1)}{\hat{i}(y_1)} \le \frac{1}{2},  $$
which contradicts to (\ref{6.42}).

\vskip 2 mm

{\bf Case 3.} {\it $\ga_1(\tau_1)$ can be connected to $N_1(1,1)\dm N_1(1,b)$ within $\Om^0(\ga_1(\tau_1))$
with $b=0$ or $1$.}

In this case, $i(y_1,1)$ and then $i(y_1)$ is even by Theorem 8.1.4 of \cite{Lon3} and Theorem 6.1. By
Theorem 8.3.1 of \cite{Lon3}, we obtain $i(y_1,m)=m(i(y_1, 1)+2)-2$ for all $m\in\N$. Thus by Theorem 6.1
we have
\bea i(y^m)=m(i(y_1)+4)-4, \qquad \forall\;m\in\N.  \nn\eea
And then $\hat{i}(y_1)=i(y_1)+4$. Because $\hat{i}(y_1)>0$ by (vii) of Claim 1, we obtain $i(y_1)\ge -2$.
By Proposition 5.4, we have $K(y_1)=1$. From (\ref{5.7}) and Remark 5.6, we then obtain
$$ \frac{\hat{\chi}(y_1)}{\hat{i}(y_1)}=\frac{k_0(y_1)-k_1(y_1)+k_2(y_1)}{i(y_1)+4}\le \frac{1}{2},  $$
which contradicts to (viii) of Claim 1.

\vskip 2 mm

{\bf Case 4.} {\it $\ga_1(\tau_1)$ can be connected to $N_1(1,1)\dm N_1(1,-1)$ within
$\Om^0(\ga_1(\tau_1))$.}

In this case, $i(y_1,1)$ and then $i(y_1)$ is odd by Theorem 8.1.4 of \cite{Lon3} and Theorem 6.1. By
Theorem 8.3.1 of \cite{Lon3}, we have $i(y,m)=m(i(y, 1)+1)-1$ for all $m\in\N$. Thus by Theorem 6.1, we
obtain
$$  i(y_1^m)=m(i(y_1)+3)-3, \qquad \forall\;m\in\N. $$
And then $\hat{i}(y_1)=i(y_1)+3$. Because $\hat{i}(y_1)>0$ by (vii) of Claim 1, we obtain
\bea  i(y_1)\ge -1.   \nn\eea
By Proposition 5.4, we have $K(y_1)=1$. Because $\nu(y_1)=2$ and Remark 5.6, we have $k_1(y_1)=1$ or 0.
Therefore we get
$$ \frac{\hat{\chi}(y_1)}{\hat{i}(y_1)}=-\frac{k_0(y_1)-k_1(y_1)}{i(y_1)+3}\le \frac{1}{2}.  $$
It contradicts to (viii) of Claim 1.

The proof of Theorem 1.4 is complete. \hfill\hb

{\bf Acknowledgment.} The authors sincerely thank the referee for his valuable comments and suggestions on this paper.

\bibliographystyle{abbrv}

\begin{thebibliography}{*}
\bibitem[BaL1]{BaL1} V. Bangert and Y. Long, The existence of two closed geodesics on every Finsler 2-sphere,
  {\it Math. Ann.} 346 (2010), 335-366.
\bibitem[BLMR]{BLMR} H. Berestycki, J. M. Lasry, G. Mancini and B. Ruf, Existence of multiple periodic orbits
  on starshaped Hamiltonian systems, {\it Comm. Pure. Appl. Math.} 38 (1985), 253-289.
\bibitem[CGH1]{CGH1} D. Cristofaro-Gardiner and M. Hutchings, From one Reeb orbit to two. arXiv:1202.4839v2,
(2012). {\it J. Diff. Geom.} to appear.
\bibitem[DDE1]{DDE1} G. Dell'Antoio, B. D'Onofrio, and I. Ekeland, Les syst\`emes hamiltoniens
convexes et pairs ne sont pas ergodiques en g\'en\'eral, {\it C. R. Acad. Sci. Paris S\'er. I}
315 (1992), 1413-1415.
\bibitem[Eke1]{Eke1} I. Ekeland, An index theory for periodic solutions of convex Hamiltonian
systems. {\it Proc. Symp. in Pure Math.} 45 (1986), 395-423.
\bibitem[Eke2]{Eke2} I. Ekeland, Convexity Methods in Hamiltonian Mechanics. Springer-Verlag.
Berlin. 1990.
\bibitem[EkH1]{EkH1} I. Ekeland and H. Hofer,  Convex Hamiltonian energy surfaces and their closed trajectories,
  {\it Comm. Math. Phys.} 113 (1987), 419-467.
\bibitem[EkL1]{EkL1} I. Ekeland and L. Lassoued,  Multiplicit\'e des trajectoires ferm\'ees d'un syst\'eme
  hamiltonien sur une hypersurface d'energie convexe. {\it Ann. IHP. Anal. non Lin\'eaire}. 4 (1987), 1-29.
\bibitem[GiG1]{GiG1} V.L. Ginzburg, Y. Goren, Iterated index and the mean Euler characteristic.
{\it Journal of Topology and Analysis}. 7 (2015), 453-481.
\bibitem[GHHM]{GHHM} V.L. Ginzburg, D. Hein, U.L. Hryniewicz, and L. Macarini, Closed Reeb orbits on the
sphere and symplectically degenerate maxima. {\it Acta Math Vietnam}. 38 (2013), 55-78.
\bibitem[Gir1]{Gir1} M. Girardi, Multiple orbits for Hamiltonian
systems on starshaped ernergy surfaces with symmetry, {\it Ann. IHP.
Analyse non lin\'eaire.} 1 (1984), 285-294.
\bibitem[GrM1]{GrM1} D. Gromoll and W. Meyer,  On differentiable functions with isolated critical points,
  {\it Topology.} 8 (1969), 361-369.
\bibitem[HWZ1]{HWZ1} H. Hofer, K. Wysocki, and E. Zehnder, The dynamics on
three-dimensional strictly convex energy surfaces. {\it Ann. of Math.} 148 (1998), 197-289.
\bibitem[HuL1]{HuL1} X. Hu and Y. Long, Closed characteristics on non-degenerate star-shaped hypersurfaces
in ${\bf R}^{2n}$, {\it Sci. China Ser. A} 45 (2002), 1038-1052.
\bibitem[LLo1]{LLo1} H. Liu and Y. Long, Resonance identity for symmetric closed characteristics on symmetric
convex  Hamiltonian energy hypersurfaces and its applications. {\it J. Diff. Equa.} 255 (2013), 2952-2980.
\bibitem[LLo2]{LLo2} H. Liu and Y. Long, The existence of two closed characteristics on every compact
star-shaped hypersurface in ${\bf R}^4$. {\it Acta Math. Sinica} (2014), DOI: 10.1007/s10114-014-4108-1.
\bibitem[LLW1]{LLW1} H. Liu, Y. Long and W. Wang, Resonance identities for closed characteristics
on compact star-shaped hypersurfaces in ${\bf R}^{2n}$. {\it J. Funct. Anal.} 166 (2014), 5598-5638.
\bibitem[LLWZ]{LLWZ} H. Liu, Y. Long, W. Wang and P. Zhang, Symmetric closed characteristics on symmetric
compact convex hypersurfaces in $\R^8$. {\it Commun. Math. Stat.} 2 (2014), 393-411.
\bibitem[LLZ1]{LLZ1} C. Liu, Y. Long and C. Zhu, Multiplicity of closed characteristics on
symmetric convex hypersurfaces in ${\bf R}^{2n}$. {\it Math. Ann.} 323 (2002), 201-215.
\bibitem[Lon1]{Lon1} Y. Long, Hyperbolic closed characteristics on compact convex smooth
hypersurfaces in ${\bf R}^{2n}$. {\it J. Diff. Equa.} 150 (1998), 227-249.
\bibitem[Lon2]{Lon2} Y. Long,  Precise iteration formulae of the Maslov-type index theory
and ellipticity of closed characteristics. {\it Advances in Math.} 154 (2000), 76-131.
\bibitem[Lon3]{Lon3} Y. Long, Index Theory for Symplectic Paths with Applications. Progress
in Math. 207, Birkh\"auser. Basel. 2002.
\bibitem[LoZ1]{LoZ1} Y. Long and C. Zhu, Closed characteristics on compact convex
hypersurfaces in ${\bf R}^{2n}$. {\it Ann. of Math.} 155 (2002), 317-368.
\bibitem[Rab1]{Rab1} P. Rabinowitz, Periodic solutions of Hamiltonian systems, {\it Comm.
Pure Appl. Math.} 31 (1978), 157-184.
\bibitem[Rad1]{Rad1}H.-B. Rademacher, Morse Theorie und geschlossene Geodatische,
{\it Bonner Math. Schriften Nr.} 229 (1992).
\bibitem[Szu1]{Szu1} A. Szulkin, Morse theory and existence of periodic solutions of convex Hamiltonian systems, {\it
  Bull. Soc. Math. France.} 116 (1988), 171-197.
\bibitem[Vit1]{Vit1} C. Viterbo, Une th\'eorie de Morse pour les syst\`emes hamiltoniens \'etoil\'es, {\it C. R. Acad.
  Sci. Paris Ser. I Math.} 301 (1985), 487-489.
\bibitem[Vit2]{Vit2} C. Viterbo, Equivariant Morse theory for starshaped Hamiltonian
systems. {\it Trans. Amer. Math. Soc.} 311 (1989), 621-655.
\bibitem[Wan1]{Wan1} W. Wang, Stability of closed characteristics on compact convex
hypersurfaces in ${\bf R}^6$. {\it J. Eur. Math. Soc.} 11 (2009), 575-596.
\bibitem[Wan2]{Wan2} W. Wang, Existence of closed characteristics on compact convex hypersurfaces in $\R^{2n}$.
arXiv:1112.5501v3, (2011).
\bibitem[Wan3]{Wan3} W. Wang, Closed trajectories on symmetric convex Hamiltonian energy
surfaces, {\it Discrete Contin. Dyn. Syst.} 32 (2012), no. 2, 679-701.
\bibitem[WHL1]{WHL1} W. Wang, X. Hu and Y. Long, Resonance identity, stability and
multiplicity of closed characteristics on compact convex hypersurfaces. {\it Duke Math.
J.} 139 (2007), 411-462.
\bibitem[Wei1]{Wei1} A. Weinstein, Periodic orbits for convex Hamiltonian systems,
{\it Ann. of Math.} 108 (1978), 507-518.
\end{thebibliography}

\end{document}